\documentclass{amsart}

\usepackage{amsmath,amsthm, amssymb,mathrsfs,  xypic, url, verbatim}
\usepackage[colorlinks]{hyperref}
\usepackage{todonotes}
\presetkeys{todonotes}{inline, size=\small}{}

 \DeclareMathOperator{\Nrd}{Nrd}

 \newcommand{\Nilp}{\mathrm{Nilp}}

\renewcommand{\AA}{\mathbb{A}}

\newcommand{\wh}[1]{\widehat{#1}}

\newcommand{\QQ}{\mathbb{Q}}
\newcommand{\ZZ}{\mathbb{Z}}
\newcommand{\FF}{\mathbb{F}}

\newcommand{\PP}{\mathbb{P}}
\newcommand{\TT}{\mathbb{T}}

\newcommand{\CC}{\mathbb{C}}

\renewcommand{\det}{\mathrm{det}}

\newcommand{\Oc}{\mathcal{O}}
\newcommand{\Ac}{\mathcal{A}}

\newcommand{\pf}{\mathfrak{p}}

\newcommand{\cont}{\mathrm{cont}}

\newcommand{\isoeq}{\cong}

\newcommand{\Tr}{\mathrm{Tr}}
\newcommand{\Gal}{\mathrm{Gal}}

\newcommand{\Aut}{\mathrm{Aut}}

\newcommand{\Frob}{\mathrm{Frob}}

\newcommand{\Spa}{\mathrm{Spa}}
\newcommand{\Spf}{\mathrm{Spf}}
\newcommand{\Spec}{\mathrm{Spec}}
\newcommand{\Hom}{\mathrm{Hom}}
\newcommand{\Isom}{\mathrm{Isom}}
\newcommand{\proet}{\mathrm{pro\acute{e}t}}
\newcommand{\Lie}{\mathrm{Lie}}

\newcommand{\GL}{\mathrm{GL}}

\newcommand{\End}{\mathrm{End}}

\newcommand{\et}{\mathrm{\acute{e}t}}

\newcommand{\ad}{\mathrm{ad}}

\newcommand{\ord}{\mathrm{ord}}
\newcommand{\Cont}{\mathrm{Cont}}

\newcommand{\Sym}{\mathrm{Sym}}

\newcommand{\LT}{\mathrm{LT}}

\newcommand{\HT}{\mathrm{HT}}

\newcommand{\abs}{\mathrm{abs}}

\newcommand{\Ig}{\mathrm{Ig}}

\newcommand{\frakm}{\mathfrak{m}}
\newcommand{\HTloc}{\mathrm{HT,loc}}
\newcommand{\ssing}{\mathrm{ss}}
\newcommand{\cusps}{\mathrm{cusps}}

\newcommand{\eval}{\mathrm{eval}}

\newcommand{\calA}{\mathcal{A}}
\newcommand{\calN}{\mathcal{N}}

\newcommand{\bbG}{\mathbb{G}}

\newcommand{\bbQ}{\mathbb{Q}}
\newcommand{\bbX}{\mathbb{X}}
\newcommand{\bbP}{\mathbb{P}}
\newcommand{\bbF}{\mathbb{F}}
\newcommand{\bbC}{\mathbb{C}}
\newcommand{\bbZ}{\mathbb{Z}}
\newcommand{\bbT}{\mathbb{T}}
\newcommand{\bbA}{\mathbb{A}}
\newcommand{\bbV}{\mathbb{V}}

\renewcommand{\l}{\left}
\renewcommand{\r}{\right}
\newcommand{\Fpbar}{\overline{\mathbb{F}}_p}
\newcommand{\mc}[1]{\mathcal{#1}}
\newcommand{\mbb}[1]{\mathbb{#1}}
\newcommand{\calI}{\mathcal{I}}
\newcommand{\calO}{\mathcal{O}}

\newcommand{\unif}{\mathrm{unif}}
\newcommand{\colim}{\mathrm{colim}}
\newcommand{\Ell}{\mathrm{Ell}}

\newcommand{\ul}[1]{\underline{#1}}

\newcommand{\Eint}{\mathfrak{E}}
\newcommand{\Ha}{\mathrm{Ha}}
\newcommand{\Yint}{\mathfrak{Y}}
\newcommand{\Xint}{\mathfrak{X}}

\numberwithin{equation}{subsection}
\theoremstyle{plain}

\newtheorem{maintheorem}{Theorem} 
 
\newtheorem{maincorollary}[maintheorem]{Corollary}

\newtheorem*{theorem*}{Theorem}

\newtheorem{theorem}[subsubsection]{Theorem}
\newtheorem{corollary}[subsubsection]{Corollary}

\newtheorem{proposition}[subsubsection]{Proposition}
\newtheorem{lemma}[subsubsection]{Lemma}

\theoremstyle{definition}

\newtheorem{example}[subsubsection]{Example}
\newtheorem{definition}[subsubsection]{Definition}
\newtheorem{remark}[subsubsection]{Remark}

\title[p-adic J-L and a question of Serre]{The spectral $p$-adic Jacquet-Langlands correspondence and a question of Serre}
\author{Sean Howe}
\address{Department of Mathematics, University of Utah, 
Salt Lake City, UT 84105}
\email{sean.howe@utah.edu}
\thanks{This material is based upon work supported by the National Science Foundation under Award
No. DMS-1704005.}

\begin{document}

\begin{abstract} We show that the completed Hecke algebra of $p$-adic modular forms is isomorphic to the completed Hecke algebra of continuous $p$-adic automorphic forms for the units of the quaternion algebra ramified at $p$ and $\infty$. This gives an affirmative answer to a question posed by Serre in a 1987 letter to Tate. The proof is geometric, and lifts a mod $p$ argument due to Serre: we evaluate modular forms by identifying a quaternionic double-coset with a fiber of the Hodge-Tate period map, and extend functions off of the double-coset using fake Hasse invariants. In particular, this gives a new proof, independent of the classical Jacquet-Langlands correspondence, that Galois representations can be attached to classical and $p$-adic quaternionic eigenforms.  
\end{abstract}

\maketitle

\setcounter{tocdepth}{2}
\tableofcontents

\newcommand{\bs}{\backslash}

\section{Introduction}
Let $p$ be a prime and let $D/\QQ$ be the (unique up to isomorphism) quaternion algebra ramified at $p$ and $\infty$. Let $\bbA$ denote the ad\`{e}les of $\bbQ$, $\bbA_f$ the finite ad\`{e}les, and $\bbA_f^{(p)}$ the finite prime-to-$p$ ad\`{e}les. Let $K^p \subset D^\times(\AA_f^{(p)})$ be a compact open subgroup. For $R$ a topological ring (e.g. $\bbC$, ${\overline{\bbF}_p}$, $\QQ_p$, or $\CC_p$), we consider the space of continuous $p$-adic automorphic forms on $D^\times$ with coefficients in $R$ and prime-to-$p$ level $K^p$, 
\[ \calA_{R}^{K^p} := \Cont(D^\times(\bbQ) \bs D^\times(\bbA) / K^p , R). \]
For $R$ totally disconnected (e.g. ${\overline{\bbF}_p}$, $\bbQ_p$, or $\bbC_p$), the archimedean component can be removed, and we have an identification 
\[ \calA_R^{K^p} = \Cont(D^\times(\bbQ) \bs D^\times(\bbA_f) / K^p , R).\]
Note that $D^\times(\bbQ) \bs D^\times(\bbA_f) / K^p$ is a profinite set. Moreover, by choosing coset representatives, it can be identified with a finite disjoint union of compact open subgroups of $D^\times(\bbQ_p)$, so that it is essentially a $p$-adic object. 

The space $\calA_{R}^{K^p}$ admits an action of the abstract double-coset Hecke algebra 
\[ \TT_{\abs} := \bbZ[ K^p \backslash D^\times(\AA_f^{(p)}) / K^p], \]
and a commuting action of $D^\times(\bbQ_p)$. In this work we, we study the spectral decomposition of $\calA_{R}^{K_p}$ under the action of $\bbT_{\abs}$. 

The classical Jacquet-Langlands correspondence \cite{jacquet-langlands:automorphic-forms-on-gl2}, proved using analytic techniques, implies that, up to twisting, the eigensystems for $\bbT_{\abs}$ acting on $\calA_{\bbC}^{K^p}$ are a strict subset of those appearing in classical complex modular forms. The eigensystem attached to a cuspidal modular form appears on the quaternionic side if and only if the corresponding automorphic representation of $\GL_2$ is discrete series at $p$.

On the other hand, arguing with the geometry of mod $p$ modular curves, Serre~\cite{serre:two-letters} showed that the eigensystems arising in $\calA_{{\overline{\bbF}_p}}^{K^p}$ are the \emph{same} as those appearing in the space of mod $p$ modular forms (see Theorem \ref{theorem:serre} below for a slight refinement of Serre's result). In particular, the gaps in the Jacquet-Langlands correspondence over $\bbC$ disappear when working mod $p$. 

The main result of this work is a natural lift of Serre's result to $\bbQ_p$: we use the geometry of the perfectoid modular curve at infinite level to show that the completed Hecke algebra of $\calA_{\bbQ_p}^{K^p}$ is isomorphic to the completed Hecke algebra of $p$-adic modular forms (see Theorem \ref{mainthm:HeckeIso} below for a precise statement). 

Theorem \ref{mainthm:HeckeIso} is compatible with the classical Jacquet-Langlands correspondence: the eigensystems appearing in classical complex quaternionic automorphic forms can be identified with the eigensystems appearing in $\calA_{\bbQ_p}^{K_p}$ such that the corresponding eigenspace contains a vector which, up to a twist, transforms via an algebraic representation of $D^\times(\bbQ_p)$ after restriction to a sufficiently small compact open subgroup. Thus, Theorem \ref{mainthm:HeckeIso} can be interpreted as saying that there is a $p$-adic Jacquet-Langlands correspondence that fills in the gaps in the classical Jacquet-Langlands correspondence. As our proof of the $p$-adic correspondence is independent of the classical correspondence, we also obtain a new proof that Galois representations can be attached to quaternionic automorphic forms (Corollary~\ref{corollary:galois-representations}). 

Both Theorems \ref{theorem:serre} and \ref{mainthm:HeckeIso} are purely \emph{spectral} Jacquet-Langlands correspondences, in the sense that they compare spectral information for a family of prime-to-$p$ Hecke operators acting on two different spaces but say little else relating the structure of these spaces; in particular,  we make no attempt here to describe the local $D^\times(\bbQ_p)$-representation appearing in a fixed Hecke eigenspace. Nevertheless, some of the methods employed in the proofs of Theorems \ref{theorem:serre} and \ref{mainthm:HeckeIso} can be used provide significant information about these local representations, and we plan to return to this in future work (see \S\ref{subsec:eigenspaces} below for further discussion).

\subsection{Serre's spectral mod $p$ Jacquet-Langlands correspondence}
Before discussing our results and techniques further, we take a detour to give a precise statement of Serre's \cite{serre:two-letters} mod $p$ correspondence.  

If we fix an isomorphism
\[ D^\times(\bbA_f^{(p)})\cong \GL_2\l(\bbA_f^{(p)}\r), \]
then we obtain an action of the Hecke algebra $\bbT_\abs$ on the space $\mc{M}_{{\overline{\bbF}_p}}^{K^p}$ of mod $p$ modular forms of prime-to-$p$ level $K^p$. In a 1987 letter to Tate, Serre \cite{serre:two-letters} proved a mod $p$ Jacquet-Langlands correspondence comparing the spectral decompositions of $\calA^{K^p}_{{\overline{\bbF}_p}}$ and $\mc{M}_{{\overline{\bbF}_p}}^{K^p}$. We state below a slight strengthening of his result, which follows essentially from Serre's proof\footnote{Actually, when $p=2$ or $3$ this method of proof leads to a small restriction on $K^p$, but this will be removed by instead deducing the mod $p$ result directly from our $p$-adic result.}. First, some notation:

Suppose $\bbT' \subset \bbT_{\abs} $ is a commutative sub-algebra and $\chi: \bbT' \rightarrow {\overline{\bbF}_p}$ is a character. Then, if $\bbT'$ acts on an ${\overline{\bbF}_p}$-vector space $V$, we may consider the $\chi$-eigenspace $V[\chi]$. If we write $\frakm_\chi := \ker \chi$, we may also consider the generalized $\chi$-eigenspace $V_{\frakm_\chi}$ (that is, the subset of elements killed by a power of $\frakm_\chi$).  

\begin{theorem}[Serre]\label{theorem:serre}
Let $\bbT' \subset \bbT_{\abs}$ be a commutative sub-algebra. Then, there is a finite collection of characters $\chi_i: \bbT' \rightarrow {\overline{\bbF}_p}$ with kernels $\frakm_i$ such that:
\begin{enumerate}
\item For each $i$, $\l(\calA^{K^p}_{{\overline{\bbF}_p}}\r)_{\frakm_i}$ and $\l(\mc{M}^{K^p}_{{\overline{\bbF}_p}}\r)_{\frakm_i}$ are non-zero; in particular, 
\[ \calA^{K^p}_{{\overline{\bbF}_p}}[\chi_i] \neq 0 \textrm{ and } \mc{M}^{K^p}_{{\overline{\bbF}_p}}[\chi_i] \neq 0, \textrm{ and } \]
\item there are direct sum decompositions
\[ \calA^{K^p}_{{\overline{\bbF}_p}} = \bigoplus_i \left( \calA^{K^p}_{{\overline{\bbF}_p}} \r)_{\frakm_{\chi_i}} \textrm{ and } \mc{M}_{{\overline{\bbF}_p}}^{K^p} = \bigoplus_i \l( \mc{M}_{{\overline{\bbF}_p}}^{K^p}\r)_{\frakm_{\chi_i}}.\]
\end{enumerate}
\end{theorem}

In other words, the Hecke eigensystems appearing in mod $p$ quaternionic automorphic forms are precisely those appearing in mod $p$ modular forms.  This stands in contrast to the classical Jacquet-Langlands correspondence, where the eigensystems appearing in quaternionic forms are a strict subset of those appearing in modular forms. The following example gives a concrete illustration:

\begin{example}\label{example:discriminant}
The discriminant form, represented by the Ramanujan series
\[ \Delta(q)=q \prod_{n \geq 1}(1-q^n)^{24} = \sum \tau(n) q^n,\]
is a weight 12 level one cuspidal eigenform whose corresponding automorphic representation is principal series at every prime $p$. Thus, the classical Jacquet-Langlands correspondence says that its Hecke eigensystem, encoded by the coefficients $\tau(\ell)$ for $\ell$ prime, does not appear in the space of classical automorphic forms on $D^\times$ for any prime $p$ (recall $p$ appears in the definition of $D^\times$).  By contrast, Theorem \ref{theorem:serre} shows that the coefficients $\tau(\ell) \mod p$ for $\ell \neq p$ are remembered by a mod $p$ quaternionic automorphic form on $D^\times$. A similar phenomenon occurs in our $p$-adic correspondence, which remembers the numbers $\tau(\ell)$ on the nose!
\end{example}

\subsection{A spectral $p$-adic Jacquet-Langlands correspondence}
\subsubsection{Serre's question} 
Serre ended his letter to Tate with a list of questions inspired by the mod $p$ Jacquet-Langlands correspondence. One of these suggests an analogous study relating $\calA_{\overline{\bbQ_p}}^{K^p}$ to $p$-adic modular forms:

\begin{quote}
\cite[paragraph (26)]{serre:two-letters} \emph{Analogues p-adiques.} Au lieu de regarder les fonctions localement constantes sur $D_\AA^\times / D_{\QQ}^\times$ \'{a} valeurs dans $\CC$, il serait plus amusant de regarder celles \`{a} valeurs dans $\overline{\QQ}_p$. Si l'on d\'{e}compose $\AA$ en $\QQ_p \times \AA'$, on leur imposerait d'\^{e}tre localement constantes par rapport \`{a} la variable dans $D_{\AA'}$ et d'\^{e}tre continues (ou analytiques, ou davantage) par rapport \`{a} la variable dans $D_p$... Y aurait-il des repr\'{e}sentations galoisiennes $p$-adiques associ\'{e}es a de telles fonctions, suppos\'{e}es fonctions propres des op\'{e}rateurs de Hecke? Peut-on interpr\'{e}ter les constructions de Hida (et Mazur) dans un tel style? Je n'en ai aucune id\'{e}e. \end{quote}

Our main result, Theorem \ref{mainthm:HeckeIso} below, shows that the answers to Serre's questions are, largely, yes. In particular, Theorem \ref{mainthm:HeckeIso} implies that Galois representations can be attached to $p$-adic quaternionic eigenforms (Corollary \ref{corollary:galois-representations} below). 

\subsubsection{A homeomorphism of completed Hecke algebras.}
The space $\calA^{K^p}_{\bbQ_p}$ of $p$-adic quaternionic automorphic forms is a $\bbQ_p$-Banach space with respect to the sup norm and the action of $\bbT_\abs$ is by bounded linear operators. For any subalgebra $\bbT'\subset \bbT_\abs$ we thus obtain a \emph{completed Hecke algebra} ${\bbT'}^\wedge_{\calA^{K^p}_{\bbQ_p}}$ by taking the closure of the image of $\bbT'$ in the algebra of bounded linear operators on $\calA^{K^p}_{\bbQ_p}$ (equipped with the topology of pointwise convergence -- see \ref{ss.completed-actions} for details). It is a topological $\bbZ_p$-algebra. 

As in the mod $p$ case, we would like to compare the Hecke action on $\calA^{K^p}_{\bbQ_p}$ to a Hecke action on a space of modular forms; in this case, we will do so by comparing completed Hecke algebras. To that end: Serre \cite{serre:formes-modulaires-et-zeta} constructed natural spaces of $p$-adic modular forms by completing spaces of classical modular forms for the $p$-adic topology on $q$-expansions (these spaces were then interpreted geometrically by Katz \cite{katz:higher-congruences, katz:l-via-moduli}). In particular, one obtains a natural $\bbQ_p$-Banach space $\mc{M}_{p-\mathrm{adic}}^{K^p}$ of $p$-adic modular forms of prime-to-$p$ level $K^p$ equipped with an action of $\bbT_\abs$ by bounded linear operators and thus a completed Hecke algebra ${\bbT'}^\wedge_{\mc{M}^{K^p}_{p-\mathrm{adic}}}$. Our main result is

\begin{maintheorem}\label{mainthm:HeckeIso}
For any sub-algebra $\TT' \subset \TT_\abs$, the identity map $\TT' \rightarrow \TT'$ extends to a canonical isomorphism of topological $\bbZ_p$-algebras
\[ {\bbT'}^\wedge_{\calA^{K^p}_{\bbQ_p}} = {\bbT'}^\wedge_{\mc{M}^{K^p}_{p-\mathrm{adic}}}. \]
\end{maintheorem}

Theorem \ref{mainthm:HeckeIso} implies Theorem \ref{theorem:serre} above (as essentially explained in \S\ref{ss.deduction-mod-p-from-p-adic} -- the point is that the maximal ideals in Theorem \ref{theorem:serre} correspond to the \emph{open} maximal ideals in the corresponding completed Hecke algebras), and gives a natural lift to characteristic zero suitable, e.g., for the construction of Galois representations. Our proof lifts Serre's approach via the mod $p$ geometry of modular curves to characteristic zero by using the $p$-adic geometry of infinite level modular curves.  

\begin{remark}\label{remark:quat-eigenform} The completed Hecke algebras do not change if we replace $\bbQ_p$ with a finite extension, or even $\bbC_p$ (and indeed this invariance under base change plays an important role in our proof). Thus, although Serre in his letter quoted above suggests the study of $\calA^{K^p}_{\overline{\bbQ}_p}$, it is natural in our setup to work over $\bbQ_p$. In particular, an eigenform in $\calA^{K^p}_{\overline{\bbQ}_p}$ will still give rise to a $\overline{\bbQ}_p$-valued character of ${\bbT'}^\wedge_{\calA^{K^p}_{\bbQ_p}}$.
\end{remark}

\subsubsection{Completed cohomology} By a result of Emerton \cite{emerton:local-global-proof} (building on work of Hida), ${\bbT'}^\wedge_{\mc{M}^{K^p}_{p-\mathrm{adic}}}$ is equal to the completed Hecke algebra of $\TT'$ acting on the completed cohomology of modular curves. On the other hand, $\Ac^{K^p}_{\QQ_p}$ \emph{is} the completed cohomology at level $K^p$ for $D^\times$. Hence we also obtain a homeomorphism between the completed Hecke algebras for the completed cohomology of $\GL_2$ and $D^\times$. In fact our proof passes first through this equivalence then uses the result of Emerton, which we establish more carefully along with some other identifications in \S\ref{ss.completed-hecke-comparison}.

\subsubsection{Galois representations}\label{sss:galois} 

Let $\TT_{\mathrm{tame}} \subset \bbT_\abs$ be the tame Hecke algebra of level $K^p$, i.e., the commutative sub-algebra generated by the Hecke operators at $\ell$ for primes $\ell\neq p$ at which $K^p$ factors as $K^{p,\ell}K_\ell$ for $K^{p,\ell} \subset D^\times(\AA_f^{(pl)})$ and $K_\ell \subset D^\times(\QQ_l)$ a maximal compact subgroup. For each such $\ell$ we write $T_\ell$ for the standard Hecke operator. Combining Theorem \ref{mainthm:HeckeIso} with known results for ${\bbT'}^\wedge_{\mc{M}^{K^p}_{p-\mathrm{adic}}}$ gives \begin{maincorollary}\label{corollary:galois-representations} 
If $\chi: \bbT^\wedge_{\mathrm{tame}, \calA^{K^p}_{\bbQ_p}} \rightarrow \bbC_p$ is a continuous character then there exists a unique semisimple continuous representation 
\[\rho: \Gal(\overline{\bbQ}/\bbQ) \rightarrow \GL_2(\bbC_p) \] 
unramified at $\ell$ as above and such that $\Tr(\rho(\Frob_\ell))=\chi(T_\ell).$ 
\end{maincorollary}

One can obtain such a $\chi$ from a quaternionic eigenform as in Remark~\ref{remark:quat-eigenform}, and thus Corollary \ref{corollary:galois-representations} attaches Galois representations to these eigenforms. 

\subsubsection{Work of Emerton}{\label{sss:emerton}
Corollary \ref{corollary:galois-representations} can also be deduced from the classical Jacquet-Langlands correspondence. In fact, a version of Theorem \ref{mainthm:HeckeIso} after localizing at a maximal ideal was first shown by Emerton \cite[3.3.2]{emerton:icm-2014} by reversing this argument: the classical correspondence gives rise to a surjective map of completed Hecke algebras
\[ {\bbT'}^\wedge_{M^{K^p}_{p-\mathrm{adic}}} \rightarrow {\bbT'}^\wedge_{\calA^{K^p}_{\bbQ_p}} \]
(which is enough to obtain Corollary \ref{corollary:galois-representations}), and then strong results in the deformation theory of Galois representations can be used to deduce that this map is an isomorphism after localizing at a maximal ideal $\frakm$ (under  minor hypotheses on $\frakm$). 

By contrast, our proof is based entirely on the $p$-adic geometry of modular curves. Thus, we obtain a new proof of Corollary \ref{corollary:galois-representations} that is independent of the classical Jacquet-Langlands correspondence, and our proof of Theorem~\ref{mainthm:HeckeIso} does not use any $R=\bbT$ theorems or other results on Galois deformations.

\subsection{Eigenspaces and the local $p$-adic Jacquet-Langlands correspondence}\label{subsec:eigenspaces}
One failing of Theorem \ref{mainthm:HeckeIso} as stated is that it says nothing about the $D^\times(\bbQ_p)$-representation appearing in the eigenspace in $\mc{A}^{K^p}_{\overline{\bbQ}_p}$ attached to a character of ${\bbT'}^\wedge_{\calA^{K^p}_{\bbQ_p}}$ valued in a finite extension of $\bbQ_p$. Indeed, because completed Hecke algebras are formed by compiling congruences of eigensystems, which may not be reflected in congruences of eigenvectors, \emph{one does not even know whether such an eigenspace is non-empty!} On the other hand, one expects that it is always non-empty, and that the $D^\times(\bbQ_p)$-representation appearing is essentially that constructed in the local correspondences of Knight \cite{knight:thesis} and Scholze \cite{scholze:lubin-tate}. 

In the course of the proof of Theorem \ref{mainthm:HeckeIso}, we produce explicit eigenvectors that show this eigenspace is non-empty at least for eigensystems attached to classical modular forms. In the author's thesis \cite{howe:thesis}, it was explained how to refine this construction so that it applies more generally to overconvergent modular forms, and it was also verified that the eigenvectors obtained are never locally algebraic for the action of $D^\times(\bbQ_p)$ (and thus in some sense are new, i.e. not obtainable by combining the classical Jacquet-Langlands correspondence and Gross's \cite{gross:algebraic-modular-forms} theory of algebraic automorphic forms, which together furnish a complete description of the locally algebraic vectors). One can now do better in at least two ways: 
\begin{enumerate} 
\item An overconvergent modular form more canonically gives rise to a non-zero map to the corresponding eigenspace in quaternionic automorphic forms from a purely local representation of $D^\times(\bbQ_p)$ constructed  as a space of distributions on the Lubin-Tate tower. The eigenvector referred to above is obtained as the image of a Dirac delta function under this map. By studying this local representation we can obtain considerably more, though still incomplete, information about the eigenspace. 
\item Combining this  with recent results of Pan \cite{pan:loc-an} on the ubiquity of overconvergent modular forms, we find that, under some minor hypotheses on the associated Galois representation, the eigenspace is always non-empty.
\end{enumerate}

Both points will be explained further in future work.

\subsection{Outline}

In \S\ref{section.preliminaries} we give some preliminaries, including in \S\ref{ss.completed-actions} some results on comparing completed Hecke algebras that will be essential in the proof of Theorem~\ref{mainthm:HeckeIso}. In \S\ref{section.modular-curves-forms-igusa} we set up our moduli problems for elliptic curves and recall the basic adelic setup for modular forms (classical, mod $p$, and $p$-adic). The main result of the section is Theorem \ref{theorem:orbit-coset}, which identifies the supersingular Igusa variety with a quaternionic coset, one of the key ingredients in the proof of both the mod $p$ and $p$-adic correspondence. Some aspects of the way we setup our moduli problems may be of independent interest; for more, we refer the reader to the introduction of \S\ref{section.modular-curves-forms-igusa}.

In \S\ref{section:mod-p} we prove the version of Serre's mod $p$ Jacquet-Langlands correspondence stated above as Theorem \ref{theorem:serre}. The proof we give is essentially that of Serre, but we emphasize carefully from the beginning the role of uniformization of the supersingular locus by the supersingular Igusa variety, which, by the above, is just a quaternionic coset. In particular, modular forms can be evaluated to quaternionic automorphic forms after choosing a trivialization of the modular bundle along the Igusa variety. A mod $p$ modular form can have zeroes or poles along the supersingular locus, but these can be cleared Hecke-equivariantly using the Hasse invariant in order to obtain a clean comparison of the corresponding Hecke algebras. When $p =2 \textrm{ or } 3$, this proof actually falls just short of the full Theorem \ref{theorem:serre}, but in \S\ref{ss.deduction-mod-p-from-p-adic} we explain how to obtain the full statement as a consequence of Theorem \ref{mainthm:HeckeIso}. 

In Section \ref{section:p-adic} we prove Theorem \ref{mainthm:HeckeIso}. Here the quaternionic coset, again in its avatar as the supersingular Igusa variety, arises naturally as a fiber of the Hodge-Tate period map in the infinite level perfectoid modular curve (as in \cite{caraiani-scholze:generic}). Thus we can evaluate modular forms to $p$-adic quaternionic automorphic forms after choosing a trivialization of the modular bundle on this fiber. A simple argument shows this evaluation map is injective; the other key property we need is density of the image. We establish this with the help of Scholze's fake Hasse invariants\footnote{We note that this step has been considerably simplified compared to the original argument in \cite{howe:thesis} by using Bhatt-Scholze's \cite{bhatt-scholze:prisms} recent result that Zariski closed equals strongly Zariski closed.}. These properties of the evaluation map are combined with the results of \S\ref{ss.completed-actions} to deduce Theorem \ref{mainthm:HeckeIso}.

\subsection{Acknowledgements}
This article grew out of the author's PhD thesis \cite{howe:thesis}, and we thank Matt Emerton, our thesis advisor, for his profound influence. We also thank Rebecca Bellovin, Ana Caraiani, Andrea Dotto, Tianqi Fan, David Hansen,
Christian Johannson, Kiran Kedlaya, Erick Knight, Daniel Le, Keerthi Madapusi-Pera, Jay
Pottharst, Peter Scholze, Joel Specter, Matthias Strauch, Jan Vonk, Jared Weinstein, and Yiwen Zhou for helpful conversations about this work. We thank an anonymous referee for helpful comments and suggestions.

\section{Preliminaries}\label{section.preliminaries}

\subsection{Numbers}

We fix a prime number $p$ and an algebraic closure $\overline{\QQ}_p$ of $\QQ_p$. We write $\bbC_p$ for the completion of $\overline{\bbQ}_p$. We denote by $\breve{\QQ}_p\subset \bbC_p$ the completion of the maximal unramified extension of $\QQ_p$ in $\overline{\QQ}_p$, and by $\breve{\ZZ}_p \subset \breve{\QQ}_p$ the completion of the ring of integers in the maximal unramified extension. We write ${\overline{\bbF}_p}$ for $\breve{\ZZ}_p/p$, an algebraically closed extension of $\FF_p = \ZZ_p/p$. There is a canonical identification $W({\overline{\bbF}_p})=\breve{\ZZ}_p$, where $W(\bullet)$ denotes Witt vectors.  

We write $\bbA$ for the ad\`{e}les of $\mbb{Q}$, $\bbA_f$ for the finite ad\`{e}les, and $\bbA_f^{(p)}$ for the finite prime-to-$p$ ad\`{e}les. We write $\hat{\bbZ}$ for the profinite completion of $\bbZ$ and $\hat{\bbZ}^{(p)}$ for the prime-to-$p$ profinite completion. We have canonical identifications 
\[ \hat{\bbZ}=\prod_{\ell \textrm{ prime} } \bbZ_\ell \textrm{ and } \hat{\bbZ}^{(p)}=\prod_{\ell \neq p \textrm{ prime }} \bbZ_\ell. \]
and inclusions $\hat{\bbZ} \subset \bbA_f$ and $\hat{\bbZ}^{(p)} \subset \bbA_f^{(p)}$ inducing isomorphisms 
\[ \hat{\bbZ} \otimes \bbQ = \bbA_f \textrm{ and } \hat{\bbZ}^{(p)} \otimes \bbQ = \hat{\bbZ}^{(p)} \otimes \bbZ_{(p)} = \bbA_f^{(p)}. \]
Here $\bbZ_{(p)}$ is interpreted via the notation $R_{\mathfrak{p}}$ for $R$ a ring and $\mathfrak{p}$ a prime ideal of $R$, which means the localization of $R$ by the multiplicative system $R\bs\mathfrak{p}.$

\subsection{Topological spaces, lisse sheaves, and torsors}\label{ss.pro-etale}

Given a topological space $T$, we denote by $\underline{T}$ the topological constant sheaf with value $T$, that is, the functor on schemes 
\[ \underline{T}(S) = \Cont(|S|, T) \]
where $|S|$ denotes the underlying topological space of $S$. When $T$ is a profinite set, it is represented by $\Spec \Cont(T, \bbZ)$, where $\bbZ$ is equipped with the discrete topology (so that the continuous functions are just locally constant). It is, in particular, a sheaf for the pro-\'{e}tale topology of \cite{bhatt-scholze:pro-etale}.

We also adopt the framework of \cite{bhatt-scholze:pro-etale} as our formalism for lisse adelic sheaves\footnote{Adelic sheaves are not discussed explicitly in loc cit., but it is no more difficult than the case of $\ell$-adic sheaves discussed in \cite[\S6.8]{bhatt-scholze:pro-etale}. Note that we are using only the most elementary parts of this formalism as we have no need for constructibility, etc. }. Thus, by a lisse $\bbA_f$-sheaf  on a scheme $S$, we mean a locally free of finite rank $\ul{\bbA_f}$ module on $S_\proet$, and similarly for $\hat{\bbZ}$, $\bbA_f^{(p)}$, $\hat{\bbZ}^{(p)}$, $\bbQ_p,$ $\bbZ_p$, etc.  A lisse $\wh{\bbZ}$-sheaf is equivalent to a compatible system of locally free $\ul{\bbZ/n\mbb{Z}}$ modules of finite rank (on either $S_\et$ or $S_\proet$ -- a locally free $\ul{\bbZ/n\mbb{Z}}$ module on $S_\proet$ is automatically classical because $\bbZ/n\bbZ$ is discrete), and similarly for $\hat{\bbZ}^{(p)}$ and $\bbZ_p$. 

For $K$ a topological group, a (right) $K$-torsor on $S_\proet$ is a sheaf $\mc{K}$ equipped with an action $\mc{K} \times \ul{K} \rightarrow \mc{K}$ such that locally on $S_\proet$, $\mc{K} \cong \ul{K}$ with the action by right multiplication. In particular, if $A=\bbA_f, \hat{\bbZ}, \bbA_f^{(p)}, \hat{\bbZ}^{(p)}, \bbZ_p,\textrm{or } \bbQ_p,$ and $V$ is a lisse $A$-sheaf of rank $n$, then 
\[ \ul{\Isom}(\ul{A}^n, V) := T \mapsto \Isom(\ul{A}^n_T, V_T) \]
 is a $\GL_n(A)$-torsor (as in these cases $\ul{\Isom}(\ul{A}^n, \ul{A}^n)= \ul{ \GL_n(A) })$. 

 The following lemma will be used as a technical tool for moving between infinite level and finite level moduli problems for elliptic curves. When $G=\GL_m(\bbZ_\ell)$ and $H=\{e\}$, it  amounts to the statement that a rank $m$ lisse $\bbZ_\ell$-sheaf is the same as a compatible family of rank $m$ lisse $\bbZ/\ell^n\bbZ$-sheafs, which will surprise nobody!
 
\begin{lemma}\label{lem.profinite-group-torsors} Let $G$ be a profinite group, $H\leq G$ a closed subgroup, and $\mc{G}$ a $G$-torsor on $\Spec R_\proet$. The map 
\begin{equation}\label{eq.extension-to-neighborhoods} \mc{H} \mapsto (\mc{H} \cdot \ul{U})_{H \leq U \leq G,\; U \textrm{compact open} } \end{equation}
is a bijection between the set of $H$-torsors in $\mc{G}$ and compatible systems of $U$-torsors in $G$ for $H \leq U \leq G,$  $U$ a compact open subgroup. 
\end{lemma}
\begin{proof}
It suffices to consider a cofinal system of $U$; thus we take a neighborhood basis of the identity in $G$ consisting of open normal subgroups $G_\epsilon$, $\epsilon \in I$, and consider only $U$ of the form $H_\epsilon := H \cdot G_\epsilon.$ Note that $\cap_{\epsilon \in I} H_\epsilon = H$. 

The key point is to show that if $(\mc{H}_\epsilon \subset \mc{G})_\epsilon$ is a compatible system of $H_\epsilon$-torsors then $\cap_{\epsilon \in I} \mc{H}_\epsilon$ admits a section on a pro-\'{e}tale cover. Indeed, then because $\cap_{\epsilon \in I}\ul{H_\epsilon} = \ul{H}$, $\cap_{\epsilon \in I} \mc{H}_\epsilon$ will automatically be an $H$-torsor, and it is straightforward to check this is a two-sided inverse to (\ref{eq.extension-to-neighborhoods}). 

For this key point, by passing to a pro-\'{e}tale cover we may assume $\mc{G}$ is trivial, i.e. we can take $\mc{G}=\ul{G}$. We now choose a compatible family of splittings of $G \rightarrow G/H_\epsilon$ (this is possible because for each $\epsilon$ the set of splittings is a finite set and the transition maps are surjective), thus we obtain a compatible family of homeomorphisms, each equivariant for the right multiplication actions of $H_\epsilon$,
\[ G = G/H_\epsilon \times H_\epsilon \]
 and, passing to the limit over $\epsilon$, a homeomorphism
\[ G = G/H \times H. \]
equivariant for the right multiplication action of $H$. 

From this it follows that the $\ul{H_\epsilon}$-torsors in $\mc{G}$ are identified with $\ul{G/H_\epsilon}(\Spec R)$, the $\ul{H}$-torsors in $\mc{G}$ are identified with $\ul{G/H}(\Spec R)$, and the map $\mc{H} \mapsto \mc{H} \cdot \ul{H_\epsilon}$ is induced by the canonical projection $G/H \rightarrow G/H_\epsilon$. The result then follows as $G/H = \lim_{\epsilon \in I} G/H_\epsilon.$ 
\end{proof}

\subsection{Elliptic curves and quasi-isogenies}
For $R$ a ring, we will consider the category $\Ell(R)$ of elliptic curves over $R$. It is $\bbZ$-linear. The isogeny category 
\[ \Ell(R) \otimes \bbQ \]
has the same objects but homomorphisms are tensored with $\bbQ$. For $E$ an elliptic curve, we sometimes write $E \otimes \bbQ$ for the corresponding element of $\Ell(R) \otimes \bbQ$, so 
\[ \Hom (E_1 \otimes \bbQ, E_2 \otimes \bbQ) = \Hom(E_1, E_2) \otimes \bbQ. \]
An \emph{isogeny} from $E_1$ to $E_2$ is a morphism $f: E_1 \rightarrow E_2$ such that $f \otimes \bbQ: E_1 \otimes \bbQ \rightarrow E_2 \otimes \bbQ$ is invertible. A \emph{quasi-isogeny} from $E_1$ to $E_2$ is an invertible morphism $f:E_1 \otimes \bbQ \rightarrow E_2 \otimes \bbQ$; we often write instead, e.g., ``$f:E_1 \rightarrow E_2$ is a quasi-isogeny." 

Similarly, we consider the prime-to-$p$ isogeny category $\Ell(R) \otimes \bbZ_{(p)}$
by replacing $\bbQ$ everywhere above with $\bbZ_{(p)}$. A \emph{prime-to-$p$ isogeny} from $E_1$ to $E_2$ is a morphism $f: E_1 \rightarrow E_2$ such that $f \otimes \bbZ_{(p)}: E_1 \otimes \bbZ_{(p)} \rightarrow E_2 \otimes \bbZ_{p}$ is invertible. A \emph{prime-to-$p$ quasi-isogeny} from $E_1$ to $E_2$ is an invertible morphism $f:E_1 \otimes \bbZ_{(p)} \rightarrow E_2 \otimes \bbZ_{(p)}$; we often write instead, e.g, ``$f:E_1\rightarrow E_2$ is a prime-to-$p$ quasi-isogeny."

\begin{remark}
When $R$ is not normal, this is not quite the category of elliptic curves up-to-isogeny (resp. prime-to-$p$ isogeny) considered in  \cite[\S3]{deligne:formes-modulaires-representations-l-adiques}, but rather the full subcategory consisting of objects with a genuine underlying elliptic curve. In general one also formally enforces effectivity of \'{e}tale descent. This full subcategory will suffice for our needs as our moduli problems typically  include rigidifying data. 
\end{remark}

\subsubsection{Tate modules} If $R/\bbQ$ (i.e. $R$ is of characteristic zero) and $E/R$ is an elliptic curve, we will consider the $p$-adic and adelic integral and rational Tate modules 
\[ T_p(E):= \lim_n E[p^n],\; V_p(E):=T_p(E)[1/p],\; T_{\hat{\bbZ}}(E) := \lim_n E[n],\;V_{\bbA_f}(E):=T_{\hat{\bbZ}}(E) \otimes \bbQ. \]
These are lisse rank two sheaves on $\Spec R$ over $\bbZ_p, \bbQ_p, \hat{\bbZ}, $ and $\bbA_f$, respectively.  All are functors on $\Ell(R)$, and $V_p$ and $V_{\bbA_f}$ factor through $\Ell(R) \otimes \bbQ$. 

If $R/\bbZ_{(p)}$ (i.e. all primes $\ell \neq p$ are invertible in $R$) and $E/\Spec R$ is an elliptic curve, then we may still form the prime-to-$p$ integral and adelic Tate modules 
\[ T_{\hat{\bbZ}^{(p)}}(E) := \lim_n E[n] \textrm{ and } V_{\bbA_f^{(p)}}(E):=T_{\hat{\bbZ}^{(p)}}(E) \otimes \bbQ = T_{\hat{\bbZ}^{(p)}}(E) \otimes \bbZ_{(p)}. \]
These are lisse rank two sheaves on $\Spec R$ over $\hat{\bbZ}^{(p)}$ and $\bbA_f^{(p)}$, respectively. Both are functors on $\Ell(R)$, and $V_{\bbA_f^{(p)}}$ factors through $\Ell(R) \otimes \bbZ_{(p)}.$

\subsubsection{Relative differentials.} For $\pi:E\rightarrow \Spec R$ an elliptic curve, $\omega_{E/R}:=\pi_* \Omega_{E/R}$ is a line bundle on $S$. Restriction induces canonical isomorphisms \[ \omega_{E/R} = 1_E^*\Omega_{E/R} = (\Lie E/R)^* \]
where here $1_E: \Spec R \rightarrow E$ is the identity section. 

The assignment $E/R \rightarrow \omega_{E/R}$ is a functor from $\Ell(R)$ to line bundles on $\Spec R$. If $R/\bbQ$, then it factors through $\Ell(R) \otimes \bbQ$, and if $R/\bbZ_{(p)}$, it factors through ${\Ell(R) \otimes \bbZ_{(p)}}$ -- indeed, for $n \in \bbZ$, the multiplication map by $n$ map $[n]:E \rightarrow E$ induces ring multiplication by $n$ on $\omega_{E/R}$, thus is invertible if  $n$ is invertible in $R$. 

\subsubsection{$p$-divisible groups}
\newcommand{\pdiv}{p-\mathrm{div}}
For $R$ a ring, a $p$-divisible group $G$ of height $h \in \ul{\mathbb{N}}(\Spec R)$ is, following Tate \cite[2.1]{tate:p-divisible-groups} and Messing \cite[I.2]{messing:crystals}, an inductive system \[ (G_i, \iota_i), i \geq 0 \]
of finite locally free commutative group schemes $G_i$ of degree $p^{ih}$ over $\Spec R$, equipped with closed immersions $\iota_i: G_i \rightarrow G_{i+1}$ identifying $G_i$ with the kernel of multiplication by $p^i$ on $G_{i+1}.$ 

We write $\pdiv(R)$ for the $\bbZ_p$-linear category of $p$-divisible groups over $R$. 
There is a natural functor 
\[ \Ell(R) \rightarrow \pdiv(R): E/R \mapsto E[p^\infty]:=(E[p^i])_i. \]
This functor factors through $\Ell(R) \otimes \bbZ_{(p)}.$ We also form the isogeny category $\pdiv(R) \otimes \bbQ_p$ and define isogenies and quasi-isogenies in the obvious way. The functor $E \mapsto E[p^\infty] \otimes \bbQ_p$ then factors through $\Ell(R) \otimes \bbQ$.

If $R$ is a $p$-adically complete ring, then we write $\Nilp_R$ for the category of $R$-algebras where $p$ is nilpotent and we view a $p$-divisible group $G=(G_i)$ over $R$ as the functor on $\Nilp_R$ 
\[ G(A) = \mathrm{colim_i} G_i(A). \]
In this case, because $E[p^\infty]_{R/p^n}$ and $E_{R/p^n}$ have the same tangent space for any $n$ and $R$ is $p$-adically complete, the functor $E \mapsto \omega_{E/S}$ factors through $E \mapsto E[p^\infty]$.

\subsection{Completing algebra actions}\label{ss.completed-actions}

In this section we develop the basic definitions for completing algebra actions and some tools for comparing completions. This material will be used in \S\ref{ss.completed-hecke-comparison}, where it is essential for the final deduction of Theorem \ref{mainthm:HeckeIso} from the key geometric input, Corollary \ref{corollary.evaluatoin-any-level}. Some results will also be used in \S\ref{ss.deduction-mod-p-from-p-adic} when we explain how Theorem \ref{theorem:serre} can be deduced from Theorem \ref{mainthm:HeckeIso}. 

The statements we give here are likely well-known to experts and the proofs are, for the most part, elementary exercises in analysis. Nonetheless  we include a full treatment because we are not aware of another suitable source in the literature.

We begin with some basic definitions in non-archimedean functional analysis. In the following, $L$ is any complete non-archimedean field.

\begin{definition} \hfill
\begin{enumerate}
\item An \emph{$L-$Banach space} is a complete topological $L$-vector space $V$ whose topology is induced by an ultrametric norm; we refer to the choice of such a norm on $V$ as a \emph{Banach norm.} 
\item A bounded collection of vectors $\{e_i \}_{i \in I}$ in an $L$-Banach space $V$ is an \emph{orthonormal basis} if every $v \in V$ can be written uniquely as
\begin{equation}\label{eq.orth-basis-vector-decomp} v = \sum_{i \in I} v_i e_i,\, v_i \in L,\, v_i \rightarrow 0.\end{equation}
\item $V$ is \emph{orthonormalizeable} if it admits an orthonormal basis. 
\end{enumerate}
\end{definition}

Note that since in (2) we assumed $\{e_i\}_{i \in I}$ was bounded, all sums of the form (\ref{eq.orth-basis-vector-decomp}) converge, and then the open mapping theorem implies that the sup norm $|v|=\sup_{i \in I} |v_i|_L$ is a Banach norm. For $L$ discretely valued (or, more generally, spherically complete), every $L$-Banach space is orthonormalizeable -- cf. \cite[Corollaire of Proposition 1 and Remarques after Proposition 2]{serre:endomorphismes-completement-continus}. 

\begin{definition}
For $V$ and $W$ two $L$-Banach spaces, we write $B(V,W)$ for the space of bounded (equivalently, continuous) linear operators from $V$ to $W$. 
\begin{enumerate}
\item The choice of Banach norms on $V$ and $W$ induces an \emph{operator norm} on $B(V,W)$ defined by 
\[ |T|=\sup_{v\in V, v\neq0} |T(v)|/|v|. \]
The operator norms for different choices of Banach norms on $V$ and $W$ are equivalent and with the induced topology $B(V,W)$ is a Banach space. 
\item The \emph{topology of pointwise convergence}\footnote{also called the \emph{strong operator topology}.} on $B(V,W)$ is defined by the family of seminorms $T \mapsto |T(v)|$ indexed by $v\in V$ (for any  Banach norm on $W$).  
\end{enumerate}
\end{definition}

The topology of pointwise convergence is uniquely determined by the property that a net $(T_j)_{j \in J}$ in $B(V,W)$ converges to $T\in B(V,W)$ if and only if $T_j(v) \rightarrow T(v)$ for all $v \in V$. It is through this characterization that we will access it. 

\begin{definition} We say $\mc{T} \subset B(V,W)$ is \emph{bounded} if it is bounded in the operator norm topology.
\end{definition} 

The following lemma is elementary but extremely useful. 
\begin{lemma}[Uniform boundedness principle]\label{lemma.uniform-boundedness} Suppose $V$ and $W$ are $L$-Banach spaces, $S \subset V$ is such that the set $L[S]$ of finite linear combinations of elements of $S$ is dense in V, $(T_{j})_{j \in J}$ is a bounded net of operators in $B(V,W)$, and $T \in B(V,W)$. Then $T_j \rightarrow T$ in the in the topology of pointwise convergence if and only if 
\begin{equation}\label{eq.net-conv-unif-bounded} \lim_{j \in J} T_j(v)=T(v) \textrm{ for all } v \in S. \end{equation}
\end{lemma}
\begin{proof}
As above, we have $T_j \rightarrow T$ in the topology of pointwise converge if and only if, for every $v \in V$, $\lim_{j \in J} T_j(v)=T(v).$ Thus one direction is immediate. For the other, suppose that (\ref{eq.net-conv-unif-bounded}) holds and fix Banach norms on $V$ and $W$. By the boundedness hypothesis, we can then choose a common bound $C\geq 1$ for the operator norms of all $T_j, j \in J$ and $T$. 

Let $v \in V$. By the density hypothesis, for any $\epsilon > 0$ we can find 
\[ v'=\ell_1 v_1 + \ldots + \ell_k v_k, v_i \in S \]
such that $|v-v'| \leq \epsilon.$ By (\ref{eq.net-conv-unif-bounded}), for each $v_i$, there is a $j_i \in J$ such that, for $j\geq j_i$, 
\[ |T(\ell_i v_i)-T_j(\ell_i v_i)|=|\ell_i||T(v_i)-T_j(v_i)| \leq \epsilon. \]
Thus, taking $j' \geq j_1, \ldots, j_k$ (the fundamental property of the directed set indexing a net is that there is always an upper bound for any finite collection of elements), we obtain that for $j \geq j'$, 
\begin{align*} |T(v)-T_j(v)| &= | (T(v') - T_j(v')) + T(v-v') - T_j(v-v')|\\  
&= \left|\sum_{i=1}^k(T(\ell_i v_i)-T_j(\ell_i v_i)) + T(v-v') - T_j(v-v')\right|\\
&\leq \max(\epsilon, |T(v-v')|, |T_j(v-v')|) \leq C\epsilon. \end{align*}
We conclude that $\lim_{j\in J} T_j(v) = T(v)$, as desired. 
\end{proof}

\begin{definition} \label{defn:StrongCompletion}
If $A$ is a ring, $L$ is a non-archimedean field, and $(W_i)_{i \in I}$ is family of $L$-Banach spaces equipped with actions of $A$ by operators in $B(W_i,W_i)$, the \emph{completion}\footnote{in the literature on Hecke algebras this is sometimes referred to as the weak completion; we avoid this terminology because of a conflict with terminology in functional analysis, where this is the completion for the strong operator topology.} of $A$ acting on $(W_i)_{i \in I}$ is the closure $\widehat{A}_{(W_i)_{i \in I}}$ of the image of $A$ in 
\[ \prod_{i \in I} B(W_i, W_i) \]
where $B(W_i, W_i)$ is equipped with the topology of pointwise convergence and the product is equipped with the product topology. Concretely, $(T_i)_{i \in I} \in \widehat{A}_{(W_i)_{i \in I}}$ if and only if there exists a net $(a_j)_{j \in J}$ in $A$ whose image converges to $(T_i)_{i \in I}$. The latter is equivalent to asking that, for any $i \in I$ and any $w_i \in W_i$, 
\[ \lim_{j \in J} a_j \cdot w_i = T_i(w_i). \]
\end{definition}

\begin{lemma}
Notation as above, if the action of $A$ on each $W_i$ is bounded (i.e. the image of $A$ in $B(W_i, W_i)$ is bounded in the operator norm topology), then $\widehat{A}_{(W_i)_{i \in I}}$	is a closed subring of $\prod_{i \in I} B(W_i, W_i)$. 
\end{lemma}
\begin{proof}
It is always a closed subgroup, so it remains just to see that under the boundedness hypothesis it is closed under composition. 

Suppose given $(T_i)_{i \in I}$ and $(S_i)_{i \in I}$ in  $\widehat{A}_{(W_i)_{i \in I}}$, and choose nets $(a_{j_T})_{j_T \in J_T}$ and $(b_{j_S})_{j_{S} \in J_S}$ whose images in $\prod_{i \in I} \End(W_i)$ converge to $(T_i)_{i \in I}$ and $(S_i)_{i \in I}$, respectively.  Then we claim that the image of
\[ (a_{j_T}b_{j_S})_{(j_T,j_S)\in J_T \times J_S} \]
converges to $(T_i \circ S_i)_{i \in I}$. It suffices to show that for any $i \in I$ and $w_i \in W_i$,
\[ \lim_{(j_T,j_S)\in J_T \times J_S} = a_{j_T}b_{j_S} \cdot w_i = T_i(S_i(w_i)). \]
We suppress the $i$s now and write $W_i=W$, $w_i=w$, $T_i=T$, $S_i=S$. 

To see the convergence, fix a Banach norm on $W$ and, by boundedness of the action, a $C\geq 1$ such that that $|a \cdot v|\leq C|v|$ for all $a \in A$ and $v \in W$. Then, 
for any $\epsilon > 0$, we may choose $j_{T,0} \in j_T$ and $j_{S,0} \in J_S$ such that 
\begin{itemize}
\item $|a_{j_T} \cdot S(w) - T(S(w))| \leq \epsilon$ for all $j_T \geq j_{T,0}$ and 
\item $|b_{j_S} \cdot w -S(w)| \leq \epsilon$ for all $j_S \geq j_{S,0}$. 
\end{itemize}
Then, for $(j_S,j_T)\geq (j_{S,0}, j_{T,0})$,
\begin{align*} |a_{j_T}b_{j_S} \cdot w - T(S(w))| &= \l|a_{j_T}\cdot(b_{j_S} \cdot w - S(w)) + (a_{j_T}\cdot S(W) - T(S(w))) \r| \\
& \leq \max(|a_{j_T}\cdot(b_{j_S} \cdot w - S(w))|,\, |a_{j_T}\cdot S(W) - T(S(w))| ) \\
& \leq \max( C\epsilon, \epsilon ) \\
& \leq C\epsilon. \end{align*}
and we conclude. 
\end{proof}

The following lemma allows for comparison to other definitions in the literature, in particular the definition given in \cite[2.1.4]{emerton:icm-2014}. 
\begin{lemma}\label{lemma.finite-dimensional-completed-Hecke}
Notation as above, suppose $L$ is discretely valued and write $\mc{O}_L$ for the ring of integers and $\pf$ for its maximal ideal. If each $W_i$ is finite dimensional and, for each $i$, $A$ preserves an $\mc{O}_L$-lattice $W_i^\circ \subset W_i$, then the action is bounded and  $\widehat{A}_{(W_i)_{i \in I}}$ is naturally identified with the closure of the image of $A$ in 
\[ \prod_{i \in I, n>0}  \End_{\mc{O}_K} (W_i^\circ / \pf^n W_i^\circ) \]
equipped with the product topology (each term is equipped with the discrete topology).
\end{lemma}
\begin{proof}
Boundedness is clear. For the rest, first note the image of $A$ in $\prod_{i \in I} B(W_i, W_i)$ factors through $\prod_{i \in I} \End_{\mc{O}_K}(W_i^\circ)$, where we identify $\End_{\mc{O}_K}(W_i^\circ)$ with the subset of $B(W_i, W_i)$ preserving $W_i^\circ$. This subset is closed, so we can form $\wh{A}_{V}$ by taking the closure of the image of $A$ in $\prod_{i \in I} \End_{\mc{O}_K}(W_i^\circ)$.  Then, for each $i$ we have 
\[ \End_{\mc{O}_K}(W_i^\circ)= \lim_n \End_{\mc{O}_K} (W_i^\circ / \pf^n W_i^\circ), \]
where each term on the right is equipped with the discrete topology. Thus 
\[ \prod_{i \in I} \End_{\mc{O}_K}(W_i^\circ) \subset \prod_{i \in I, n>0}  \End_{\mc{O}_K} (W_i^\circ / \pi^n) \]
is closed so we may compute $\wh{A}_V$ by taking the closure in the space on the right. 
\end{proof}

The following lemma says that completion is insensitive to base extension. This is useful as our comparisons of Hecke-modules take place over very large extensions of $\QQ_p$, whereas one is typically interested in Hecke algebras over $\ZZ_p$.

\begin{lemma}\label{lem:CompBaseChange}
Let $L \subset L'$ be an extension of complete non-archimedean fields, and let $A$ be a ring. Suppose $(W_i)$ is a family of orthonormalizable $L$-Banach spaces equipped with bounded actions of $A$. Then the identity map $A \rightarrow A$ extends uniquely to a topological isomorphism 
\[ \widehat{A}_{(W_i)_{i\in I}} = \widehat{A}_{(W_i \widehat{\otimes}_L L')_{i \in I}}. \]
\end{lemma}  
\begin{proof}
Immediate by applying the uniform boundedness principle (Lemma \ref{lemma.uniform-boundedness}) to an orthonormal basis and using the fact that an orthonormal basis remains an orthonormal basis under completed base change. 
\end{proof}

The following is our main technical tool for comparing completed Hecke algebras. 

\begin{lemma}\label{lem:CompDense}
Suppose $V$ is an orthonormalizable $L$-Banach space equipped with a bounded action of a ring $A$, and $(W_i)_{i \in I}$ is a collection of $A$-invariant closed subspaces such that the span of $\bigcup_i W_i$ is dense in $V$. Then $\widehat{A}_{(W_i)_{i\in I}} = \widehat{A}_{V}$.

\end{lemma} 

\begin{remark}
In this setup, each $W_i$ is automatically a Banach space as a closed subspace of a Banach space and the action on $W_i$ is automatically bounded.
\end{remark}

\begin{proof}
We abbreviate $\widehat{A}_W := \widehat{A}_{(W_i)_{i\in I}}\subset \prod B(W_i, W_i) $. We then obtain a map $\widehat{A}_V \rightarrow \widehat{A}_W$ via restriction: if the image of $(a_j)_{j \in J}$ in $B(V,V)$ converges to $T$, then in particular the image of $(a_j)_{j \in J}$ in $B(W_i, W_i)$ converges to $T|_{W_i}$. This restriction map is injective by the density hypothesis.

We show now that it is surjective. The key observation that makes this possible is that, by the density hypothesis, we may choose an orthonormal basis $(e_m)_{m \in M}$ for $V$ consisting of vectors $e_m$ each of which is a finite linear combination of vectors in the subspaces $W_i$: indeed, if we fix a pseudo-uniformizer $\pi$ in $\mc{O}_L$ and an arbitrary orthonormal basis $(f_m)_{m \in M}$, then any collection of vectors $(e_m)_{m \in M}$ with $|f_m - e_m| \leq |\pi|$ will also be an orthonormal basis. 

Now, suppose $(T_i)_{i \in I} \in \widehat{A}_W$, and fix a net $a_j$ of elements of $A$ whose image converges to $(T_i)_{i \in I}$. Then we find that for each $m$, $\lim_j a_j \cdot e_m$ exists in $V$, call it $v_m$ and, by boundedness of the action, the set of $v_m$ is bounded. There is thus a unique bounded linear operator $T:V\rightarrow V$ such that $T(e_m)=v_m$. We conclude by the uniform boundedness principle that the image of $a_j$ in $B(V,V)$ converges to $T$, and then by restriction that in fact $T|_{W_i}=T_i$. 

Now, a topology is uniquely determined by the knowledge of which nets converge to which points. With this bijection established, the uniform boundedness principle tells us that the same nets converge to the same points, so the bijection is a homeomorphism. 
\end{proof}

The following lemma combines some of the above results, and will be used in \S\ref{ss.deduction-mod-p-from-p-adic} to deduce Theorem \ref{theorem:serre} from Theorem \ref{mainthm:HeckeIso}. 
\begin{lemma}\label{lemma:dense-prodiscrete}
Suppose $L$ is discretely valued and write  $
\mc{O}_L$ for the ring of integers and $\pf$ for the maximal ideal. Suppose $V$ is an $L$-Banach space and $A$ acts on $V$ preserving a bounded open $\mc{O}_K$-lattice $V^\circ$, and $(W_i)_{i \in I}$ is a filtered system of finite dimensional $A$-invariant subspaces such that $\bigcup_{i \in I} W_i$ is dense in $V$. Then, writing 
\[ W_i^\circ = W_i \cap V^\circ, \;\overline{W}_{i,n}=W_i ^\circ / \pf^n W_i^\circ = W_i^\circ/ W_i \cap \pf^n V^\circ, \]
and $A_{i,n}$ for the image of $A$ in $\End_{\mc{O}_K}( \overline{W}_{i,n})$, we have
\[ \wh{A}_V = \wh{A}_{(W_i)_{i \in I}} = \lim_{(i,n) \in I \times \mbb{N}} A_{i,n} \]
 where each term in the limit is equipped with the discrete topology. 
\end{lemma} 
\begin{proof}
It follows from Lemma \ref{lem:CompDense} that $\wh{A}_V = \wh{A}_{(W_i)_{i \in I}}$. The result then follows from Lemma \ref{lemma.finite-dimensional-completed-Hecke}, since in this case the closure of the image of $A$ will be identified with the limit of the $A_{i,n}$ as a subset of the product appearing there. 
\end{proof}

\section{Modular curves and Igusa varieties}\label{section.modular-curves-forms-igusa}

In this section we study some moduli problems for elliptic curves. In \S\ref{ss.modular-curves} we give isogeny formulations for some classical moduli problems and recall the standard representability results. In \S\ref{ss.modular-forms} we recall the construction of the modular bundle and the adelic representations on modular forms, as well as the construction of the Hasse invariant. In \S\ref{ss.supersingular-and-ordinary-loci} we recall the construction of the supersingular and ordinary loci on the mod $p$ modular curve. In \S\ref{ss.ordinary-igusa-mod-p-padic} we recall some Igusa moduli problems over the ordinary locus and their relation with mod $p$ and $p$-adic modular forms as developed by Katz \cite{katz:higher-congruences}. 

In \S\ref{ss.igusa-variety}-\S\ref{ss.quat-coset} we undertake a study of the supersingular Igusa variety, culminating with the identification of the supersingular Igusa variety with a quaternionic coset in Theorem \ref{theorem:orbit-coset}. Everything here except this final identification is a very special case of results of Caraiani-Scholze \cite{caraiani-scholze:generic}. However, following our treatment of modular curves, we take a resolutely ``top-down" approach, and for the most part\footnote{We do however appeal crucially in Lemma \ref{lemma.iso-isog-p-div} to the construction of an internal hom for $p$-divisible groups over $\Fpbar$ as established in \cite{caraiani-scholze:generic} based on work of Chai-Oort.} our treatment here is independent of the results of \cite{caraiani-scholze:generic}. We lean instead on the Hasse invariant and other ideas specific to this special case. 

As to the identification with a quaternionic coset, the basic idea is already present in Serre's letter \cite{serre:two-letters}, so our main contribution is a careful treatment by exploiting the group action at infinite level. This identification is a key ingredient in both the mod $p$ correspondence in \S\ref{section:mod-p} and the $p$-adic correspondence in \S\ref{section:p-adic}.

Finally we remark that, motivated by our specific needs, we have made what appear to be some nonstandard choices in defining our moduli problems:
\begin{enumerate}
\item We allow level defined by an arbitrary \emph{closed} adelic subgroup, which facilitates the free usage of large group actions on infinite level moduli problems and, in particular, transparent passage between the infinite level prime-to-$p$ moduli problem over $\bbZ_{(p)}$ and infinite level moduli problem over $\bbQ$. 
\item We give an up to isogeny definition of level structure that does not require the base scheme to be locally noetherian (i.e. does not use the (pro)-\'{e}tale fundamental group). In particular, this is necessary to allow arbitrary closed subgroups as above, but also allows us to evaluate on, e.g., perfectoid rings and other very non-noetherian objects without appealing behind the scenes to noetherian approximation.  
\end{enumerate}
We accomplish both of these goals by interpreting the sentence ``level $K$ structure on $E$ is a $K$-orbit of trivializations of $V_{\bbA_f}(E)$" literally, i.e. as the choice of a $K$-torsor in $\ul{\Isom}(\ul{\bbA_f}^2, V_{\bbA_f}(E))$. All representability statements are deduced from  classical results on finite level curves, and ultimately all of our arguments could be run in a more classical setup, as the diligent reader will have no trouble verifying.

\subsection{Modular curves}\label{ss.modular-curves}

\begin{definition}[The level $K$ elliptic moduli functor] Let $K \subset \GL_2(\bbA_f)$ be a closed subgroup. Let $Y_K$ be the functor on $\bbQ$-algebras
\[ Y_K: R \mapsto \{ (E, \mc{K}) \} /\sim \]
sending $R/\bbQ$ to the set of equivalence classes of pairs $(E, \mc{K})$ where
\begin{enumerate}
\item $E/R$ is an elliptic curve,
\item $\mc{K} \subset \ul{\Isom}( (\ul{\bbA_f})^2, V_{\bbA_f }(E))$ is a $K$-torsor,
\item The relation $\sim$ is defined by $(E,\mc{K}) \sim (E', \mc{K}')$ if there is a quasi-isogeny $q: E \rightarrow E' $ such that 
$q(\mc{K}) = \mc{K}'.$
\end{enumerate}
The topological constant sheaf on the normalizer of $K$, $\ul{N_{\GL_2(\bbA_f)}(K)}$, acts on $Y_K$, and for $K_1 \leq K_2$ we have the obvious map 
\[ Y_{K_1} \rightarrow Y_{K_2},  \; (E, \mc{K}_1) \mapsto (E, \mc{K}_1 \cdot \ul{K_2} ). \]
\end{definition}

\begin{example}[Infinite level]\label{example.infinite-level} Take $K=\{e\}$. Then the $K$-torsor $\mc{K}$ appearing in the moduli problem $Y_{K}$ is just a section of 
	\[ \ul{\Isom}( \ul{\bbA_f}^2, V_{\bbA_f}(E)), \]
	i.e. an isomorphism $\varphi_{\bbA_f}:\ul{\bbA_f}^2 \xrightarrow{\sim} V_{\bbA_f}(E)$, and the condition in the equivalence relation becomes $q \circ \varphi_{\bbA_f} = \varphi'_{\bbA_f}.$ The group action is by all of $\ul{\GL_2(\bbA_f)}$, and in this notation it acts by composition with $\varphi_{\bbA_f}.$ \textbf{When $K=\{e\}$, we will typically omit it from the notation and write simply $Y=Y_{\{e\}}$. }
\end{example}

Removing any level structure at $p$, we obtain a variant over $\bbZ_{(p)}$:

\begin{definition}[The integral level $K^p$ elliptic moduli functor] Let $K^p \leq \GL_2(\bbA_f^{(p)})$ be a closed subgroup. Let $\Yint_{K^p}$ be the functor on $\bbZ_{(p)}$-algebras 
\[ \Yint_{K^p}: R \mapsto \{ (E, \mc{K}^p) \} /\sim \]
sending $R/\mbb{Z}_{(p)}$ to the set of equivalence classes of pairs $(E, \mc{K}^p)$ where
\begin{enumerate}
\item $E/R$ is an elliptic curve,
\item $\mc{K}^p \subset \ul{\Isom}( (\ul{\bbA_f^{(p)}})^2, V_{\bbA_f^{(p)}}(E))$ is a $K^p$-torsor,
\item The relation $\sim$ is defined by $(E, \mc{K}^p) \sim (E',  {\mc{K}^p}')$ if there is a prime-to-$p$ quasi-isogeny $q: E \rightarrow E' $ such that  $q(\mc{K}^p) = {\mc{K}^p}'.$
\end{enumerate}
The topological constant sheaf on the normalizer of $K$, $\ul{N_{\GL_2(\bbA_f^{(p)})}(K^p)}$, acts on $\Yint_{K^p}$, and for $K^p_1 \leq K^p_2$ we have the obvious map 
\[ \Yint_{K^p_1} \rightarrow \Yint_{K^p_2},  \; (E, \mc{K}_1^p) \mapsto (E,  \mc{K}_1^p \cdot \ul{K^p_2}). \]
\end{definition}

\begin{example}[Integral infinite level]\label{example.integral-infinite-level}
As in Example \ref{example.infinite-level}, when $K^p=\{e\}$, $\mc{K}^p$ is simply the choice of an isomorphism $\varphi_{\bbA_f^{(p)}}:\ul{\bbA_f^{(p)}}^2 \xrightarrow{\sim} V_{\bbA_f^{(p)}}(E)$. \textbf{When $K^p=\{e\}$, we will typically omit it from the notation and write simply $\Yint=\Yint_{\{e\}}.$}  
\end{example}

Arguing as in \cite[Corollaire 3.5]{deligne:formes-modulaires-representations-l-adiques}, we find

\begin{lemma}\label{lemma:generic-integral-fiber} Let $K^p \leq \GL_2(\bbA_f^{(p)})$ be a closed subgroup. The assignment
\[ (E, \mc{K}^p) \rightarrow (E, \ul{\Isom}(\bbZ_p^2, T_p E) \times \mc{K}^p) \]
induces an isomorphism 
\[ \Yint_{K^p, \bbQ} \xrightarrow{\sim} Y_{\GL_2(\bbZ_p)K^p}. \]
\end{lemma} 

\begin{example} Lemma \ref{lemma:generic-integral-fiber} gives $\Yint_{\bbQ}=Y_{\GL_2(\bbZ_p)},$
where on the right $\GL_2(\bbZ_p)$ is viewed as a closed subgroup of $\GL_2(\bbA_f)$, and, as above, $\Yint=\Yint_{\{e\}}$. This identification explains one reason why it is convenient to allow an arbitrary closed subgroup in the formulation of the moduli problem.  
\end{example}

\begin{definition}
A closed subgroup $K \leq \GL_2(\bbA_f)$ (resp. $K^p \leq \GL_2(\bbA_f^{(p)})$) is \emph{sufficiently small} if it stabilizes a $\hat{\bbZ}$-lattice $\mc{L} \subset \bbA_f^2$ (resp. a $\hat{\bbZ}^{(p)}$-lattice $\mc{L} \subset (\bbA_f^{(p)})^2$) and lies in the kernel of the map $\GL(\mc{L}) \rightarrow \GL(\mc{L}/n\mc{L})$ for some $n \geq 3$ (resp. and $(n,p)=1$).  
\end{definition}

Note that if $K_2 \leq K_1 \leq \GL_2(\bbA_f)$ (resp. $K_2^p \leq K_1^p \leq \GL_2(\bbA_f^{(p)})$),  are closed subgroups and $K_1$ (resp. $K_1^p$) is sufficiently small then so is $K_2$ (resp. $K_2^p$). Moreover, the property of being a sufficiently small closed subgroup of $\GL_2(\bbA_f)$ (resp. $\GL_2(\bbA_f^{(p)})$) is preserved under conjugation by $\GL_2(\bbA_f)$ (resp. $\GL_2(\bbA_f^{(p)})$). Because any lattice $\mc{L}$ is in the $\GL_2(\bbA_f)$-orbit of $\hat{\bbZ}^2$ (resp. $\GL_2(\bbA_f^{(p)})$-orbit of $(\hat{\bbZ}^{(p)})^2$), being sufficiently small is equivalent to being contained in a conjugate of the standard principal congruence subgroup of level $n \geq 3$ (resp. and $(n,p)=1$). 

The main representability results are:
\begin{proposition}\label{prop:modular-curves} If $K \leq \GL_2(\bbA_f)$ (resp. $K^p \leq \GL_2(\bbA_f^{(p)})$) is a sufficiently small closed subgroup then $Y_K$ (resp. $\Yint_{K^p}$) is represented by an affine scheme over $\Spec \bbQ$ (resp. $\Spec \bbZ_{(p)}$), and the natural map 
\begin{equation}\label{eq.map-to-limit} Y_K \rightarrow \lim_{\substack{K' \textrm{ compact open}\\ K \leq K' \leq \GL_2(\bbA_f)}} Y_{K'} \;\; \textrm{ (resp. } \Yint_{K^p} \rightarrow \lim_{\substack{{K^p}' \textrm{ compact open}\\ {K^p} \leq {K^p}' \leq \GL_2(\bbA_f^{(p)})}} \Yint_{{K^p}'} \textrm{)} \end{equation}
is a $\ul{N_{\GL_2(\bbA_f)}(K)}$-equivariant (resp. $\ul{N_{\GL_2(\bbA_f^{(p)})}(K^p)}$-equivariant) isomorphism, where the action on the right-hand side is induced by the action on the tower that permutes the terms by conjugation (i.e. right multiplication by $h$ sends $Y_{K'}$ to $Y_{h^{-1}K'h}$). 

If $K$ (resp. $K^p$) is furthermore compact open then $Y_K$ (resp. $\Yint_{K^p}$) is a smooth affine curve. Moreover, for $K_1 \leq K_2$ (resp $K_1^p \leq K_2^p$) sufficiently small closed subgroups the natural map $Y_{K_1} \rightarrow Y_{K_2}$ (resp. $\Yint_{K_1^p} \rightarrow \Yint_{K_2^p}$) is profinite \'{e}tale, and, if $K_1 \trianglelefteq K_2$ (resp. $K_1^p \trianglelefteq K_2^p$) it is Galois with group $K_1/K_2$ (resp. $K_1^p/K_2^p$). 
\end{proposition}
\begin{proof}
We argue only in the case over $\bbQ$, as the argument over $\bbZ_{(p)}$ is essentially the same. If we fix a $\hat{\bbZ}$ lattice $\mc{L}\subset \bbA_f^{2}$ preserved by $K$ then, as in \cite[Corollaire 3.5]{deligne:formes-modulaires-representations-l-adiques}, we see that the moduli problem can be replaced with an equivalent up to isomorphism moduli problem by taking $\mc{K}$ in $\ul{\Isom}(\mc{L}, T_{\hat{\bbZ}^{\bullet}}(E))$. The assertion that (\ref{eq.map-to-limit}) is an isomorphism then amounts to the following: for $E/R$ an elliptic curve, if we consider the $\ul{\GL(\mc{L})}$-torsor  $\mc{G}:=\ul{\Isom}(\mc{L}, T_{\hat{\bbZ}^{\bullet}}(E))$, we must show that the following data are equivalent:
\begin{enumerate}
\item a $K$--torsor inside  $\mc{G}$,
\item a system of $K'$-torsors inside $\mc{G}$ for $K'$ compact open, $K \leq K' \leq \GL(\mc{L})$, compatible under inclusion. 
\end{enumerate}
This equivalence is provided by Lemma \ref{lem.profinite-group-torsors}.

Since we have established that (\ref{eq.map-to-limit}) is an isomorphism, the rest of the claim for general $K$ is essentially formal if we can establish the representability claims for $K$ compact open. But for $K$ compact open we can conjugate to assume the lattice $\mc{L}$ as above is $\hat{\bbZ}^2$, and then the representability statements  are consequences of the classical theory of finite level modular curves as in, e.g. \cite{katz-mazur:moduli}. 
\end{proof}

\begin{definition}[Compactified modular curves] 
For $K \leq \GL_2(\bbA_f)$  (resp. $K^p \leq \GL_2(\bbA_f^{(p)})$) a sufficiently small compact open subgroup, we form  compactifications $X_K$ (resp. $\Xint_{K^p}$) as in \cite[8.6]{katz-mazur:moduli} after fixing a lattice $\mc{L} \subset \bbA_f^2$ (resp. $\mc{L} \subset (\bbA_f^{(p)})^2$) preserved by $K$ (resp. $K^p$) to relate to classical finite level moduli problems as in the proof above. We obtain smooth projective curves $X_K/\bbQ$ (resp. $\Xint_{K^p}/\bbZ_{(p)}$), and the finite \'{e}tale maps in the tower of $Y_K$ (resp. $\Yint_{K^p}$) for $K$ (resp. $K^p$) sufficiently small compact open extend to finite maps in the tower of $X_K$ (resp. $\Xint_{K^p}$). 

It can be checked that the natural group actions also extend. For $X_K$ and $\Xint_{K^p,\bbF_p}$ this is even immediate because the smooth compactifications of smooth curves are functorial over a perfect field -- using this we could also \emph{define} $X_K$ and $\Xint_{K^p,\bbF_p}$  with no reference to the moduli problem. 

We extend these definitions to $K$ (resp. $K^p$) sufficiently small closed by taking limits as in Proposition \ref{prop:modular-curves}, and the resulting objects are schemes because the transition maps are affine. 

\end{definition}

We refer to the boundary $X_K\backslash Y_K$ (resp. $\Xint_{K^p} \backslash \Yint_{K^{p}}$) with its reduced subscheme structure as the cusps. The cusps can also be described as filling in the punctures corresponding to level structure on the Tate curve as in \cite[8.11]{katz-mazur:moduli}.

\subsection{Modular forms}\label{ss.modular-forms}
For any sufficiently small closed $K \leq \GL_2(\bbA_f)$ (resp. $K^p \leq \GL_2(\bbA_f^{(p)})$), we have a universal elliptic curve $E_K/Y_K$ (resp. $\Eint_{K^p} / \Yint_{K^p}$), determined up to unique isogeny (resp. prime-to-$p$ isogeny). We write simply $\omega$ for the line bundle $\omega_{E/K} / Y_K$ (resp. $\omega_{\Eint_{K^p}/\Yint_{K^p}}/\Yint_{K^p}$) -- it is determined up to unique isomorphism compatibly with all pullbacks, base change, and group actions discussed so far, so that this notation will cause no confusion.  

For every sufficient small compact open $K$ (resp. $K^p$), we extend $\omega$ to $X_K$ (resp. to $\Xint_{K^p}$) in the standard way by allowing sections with holomorphic $q$-expansions at each cusp. Direct computation shows this is compatible with all pullbacks, base change, and group actions discussed so far, so that we can extend this definition to any sufficiently small closed $K$ (resp. $K^p$) and again no confusion will be caused by referring to the extended line bundle also as $\omega$. 

We consider the smooth $\GL_2(\bbA_f)$-representation of \emph{modular forms}, 
\[ M_{k,\bbQ} = H^0(X, \omega^k). \]
For $K$ any sufficiently small closed subgroup, pullback from level $K$ identifies 
\[ M_{k,\bbQ}^K = H^0(X_K, \omega^k). \]
Applied to $K$ sufficiently small compact open, we deduce that $M_{k,\bbQ}$ is an admissible representation of $\GL_2(\bbA_f)$. Applied to $K=\GL_2(\bbZ_p)$, we obtain
\[ M_{k,\bbQ}^{\GL_2(\bbZ_p)} = H^0(X_{\GL_2(\bbZ_p)}, \omega^k)=H^0(\Xint_{\bbQ}, \omega^k)=H^0(\Xint, \omega^k)\otimes_{\bbZ_{(p)}} \bbQ. \]
In particular, if we write 
\[ M_{k,\bbF_p} = H^0(\Xint_{\bbF_p}, \omega^k), \]
an admissible $\bbF_p$-representation of $\GL_2(\bbA_f^{(p)})$ by the same argument as above, then $H^0(\Xint, \omega^k)$ is a natural $\bbZ_{(p)}$-lattice in $M_{k,\bbQ}^{\GL_2(\bbZ_p)}$ equipped with a $\GL_2(\bbA_f^{(p)})$-equivariant reduction map to $M_{k,\bbF_p}$.
\begin{remark} Neither $M_{k,\bbF_p}$ nor the image of reduction is what is typically referred to as mod $p$ modular forms! We will recall this definition in \S\ref{ss.ordinary-igusa-mod-p-padic} below. 
\end{remark}

\subsubsection{The Hasse invariant}\label{sss.hasse} We now recall how, to any elliptic curve $E/R$ for $R/\bbF_p$, one can attach a canonical section $\Ha(E/R) \in \omega^{p-1}_{E/R}$, the Hasse invariant. We follow one of the approaches described in \cite[12.3]{katz-mazur:moduli}.

Since the section $\Ha(E/R)$ can be constructed Zariski locally, it suffices to assigns to any pair $(E/R, \alpha)$ where $R$ is an $\mbb{F}_p$-algebra, $E/ R$ is an elliptic curve and $\alpha \in \omega_{E/R}$ is a nonvanishing invariant differential, an element $\Ha(E/ R,\alpha)$ of $R$ such that, for $a \in A^\times$,  
\[ \Ha(E/ R, a \alpha)= a^{-(p-1)}\Ha(E/ R, \alpha) \]
and whose formation is functorial in base change and isomorphism. In this case, to give our rule we first take the invariant derivation $\partial_{\alpha}$ that is dual to $\alpha$, then form 
\[ \partial_{\alpha}^p := \underbrace {\partial_{\alpha} \circ \ldots \circ \partial_{\alpha} }_{\textrm{$p$ times}}, \]
which is also an invariant derivation and thus a multiple of $\partial$. Then the equation 
\[ \partial^p_\alpha = \Ha(E/ R, \alpha) \partial_\alpha, \]
defines $\Ha(E/ R, \alpha)$, and it is straightforward to check this satisfies the desired transformation rule if we scale $\alpha$ and is functorial in base change and isomorphism. 

In fact, the construction is also functorial in prime-to-$p$ quasi-isogenies: it suffices to observe that a prime-to-$p$ quasi-isogeny induces an isomorphism of $p$-divisible groups and in particular of formal groups, and that the action of an invariant derivation is completely determined by its action on the formal group. This observation also has the important consequence that the resulting section $\Ha(E/R)$ can be constructed entirely in terms of $E[p^\infty]$. 

Applying this construction to the universal elliptic curve over $\Yint_{\bbF_p}$, we obtain a $\GL_2(\bbA_f^{(p)})$-invariant section of $\omega^{p-1}$. A direct computation on the Tate curve (see \cite[Theorem 12.4.2]{katz-mazur:moduli}) shows that its $q$-expansions are constant equal to $1$ at every cusp, thus it extends to \[ \Ha \in M_{p-1, \bbF_p}^{\GL_2(\bbA_f^{(p)})}. \]

\subsection{The supersingular and ordinary loci}\label{ss.supersingular-and-ordinary-loci}

Let 
\[ b_{\ssing}= \begin{pmatrix} 0 & p \\ 1 & 0 \end{pmatrix} \in M_2(\bbZ_p) \textrm{ and } b_\ord=\begin{pmatrix} p & 0 \\ 0 & 1 \end{pmatrix} \in M_2(\bbZ_p).\]
Let $\bbX_{\ssing} / \bbF_p$ (resp. $\bbX_{\ord}/\bbF_p$) be the $p$-divisible group corresponding to the covariant Dieudonn\'{e} module $\bbZ_p^2$ with Frobenius $F$ acting by $b_{\ssing}$ (resp. $b_\ord$). Then $\bbX_{\ssing}$ is a connected one-dimensional height two $p$-divisible group, while $\bbX_{\ord}=\mu_{p^\infty} \times \bbQ_p/\bbZ_p$ is a one-dimensional height two $p$-divisible group with non-trivial \'{e}tale part. It follows from the classification of $p$-divisible groups by Dieudonn\'{e} modules that, for any algebraically closed $\kappa/\bbF_p$, every height two one-dimensional $p$-divisible group over $\kappa$ is quasi-isogenous/isomorphic\footnote{The fact that quasi-isogeny classes are equal to isomorphism classes in this case is, of course, one of the many facts that makes $\GL_2$ very special!} to exactly one of $\bbX_{\ord, \kappa}$ and $\bbX_{\ssing, \kappa}$. 

In particular, this applies to $E[p^\infty]$ for $E/\kappa$ an elliptic curve. It thus makes sense, for any $K^p$ sufficiently small, to define the supersingular (resp. ordinary) locus $\Yint_{K^p,\bbF_p}^\ssing$ (resp. $\Yint_{K^p,\bbF_p}^\ord$) in $\Yint_{K^p,\bbF_p}$ as the locus whose geometric points are such that the $p$-divisible group of universal elliptic curve is quasi-isogenous/isomorphic to $\bbX_\ssing$ (resp. $\bbX_\ord$). We then have 
\[ \Yint_{K^p, \bbF_p} = \Yint_{K^p,\bbF_p}^\ord \sqcup \Yint_{K^p,\bbF_p}^\ssing \]
and it can be shown that the supersingular locus is closed (and thus the ordinary locus is open). In fact, a local computation as in \cite[Theorem 12.4.3]{katz-mazur:moduli} gives 
\begin{lemma}\label{lemma.Hasse-vanishing} For $K^p$ sufficiently small, the vanishing locus of $\Ha$, viewed as section of $H^0(\Xint_{K^p}, \bbF_p)$, is $\Yint_{K^p,\bbF_p}^\ssing$ with its reduced closed subscheme structure.
\end{lemma}

We write also
\[ \Xint_{K^p,\bbF_p}^\ssing = \Yint_{K^p,\bbF_p}^\ssing = V(\Ha) \textrm{ and } \Xint_{K^p,\bbF_p}^\ord=\Xint_{K^p,\bbF_p} \backslash \Xint_{K^p,\bbF_p}^\ssing. \]

\subsection{Ordinary Igusa varieties, mod $p$, and $p$-adic modular forms}\label{ss.ordinary-igusa-mod-p-padic}
Following Katz \cite{katz:higher-congruences}, we will \emph{define} mod $p$ and $p$-adic modular forms as functions on certain moduli problems. This has the advantage of rendering group actions transparent, and mirrors our approach to the supersingular Igusa variety and quaternionic $p$-adic automorphic forms. Functions on these moduli problems have $q$-expansions, and the connection with the completion of $q$-expansions of classical modular forms as in Serre's \cite{serre:formes-modulaires-et-zeta} perspective then comes from evaluation along a canonical trivialization of $\omega$ and the $q$-expansion principle. 

We begin by treating mod $p$ modular forms.  
\begin{definition}[The $\mu_p$-Igusa moduli problem] Let $K^p \leq \GL_2(\bbA_f^{(p)})$ be a closed subgroup. Let $\Ig^{\mu_p}_{K^p}$ be the functor on $\bbF_{p}$-algebras 
\[ \Ig^{\ord}_{K^p,\mu_p}: R \mapsto \{ (E, \varphi_p, \mc{K}^p) \} /\sim \]
sending $R/\bbF_p$ to the set of equivalence classes of triples $(E, \varphi_p, \mc{K}^p)$ where
\begin{enumerate}
\item $E/R$ is an elliptic curve,
\item $\varphi_p: \mu_p \xrightarrow{\sim} \hat{E}[p]$ is an isomorphism,
\item $\mc{K}^p \subset \ul{\Isom}( (\ul{\bbA_f^{(p)}})^2, V_{\bbA_f^{(p)}}(E))$ is a $K^p$-torsor,
\item The relation $\sim$ is defined by $(E, \mc{K}^p) \sim (E',  {\mc{K}^p}')$ if there is a prime-to-$p$ quasi-isogeny $q: E \rightarrow E' $ such that 
$q\circ \varphi_p = \varphi_p'$ and $q(\mc{K}^p) = {\mc{K}^p}'.$
\end{enumerate}
\end{definition}

\textbf{As usual, when $K^p=\{e\}$ we drop it from the notation.}
For $K^p$ sufficiently small, the functor $\Ig^{\ord}_{K^p,\mu_p}$ is a finite \'{e}tale $(\bbZ/p\bbZ)^\times$ cover of $\Yint^{\ord}_{K^p, \bbF_p}$, where $(\bbZ/p\bbZ)^\times$ acts by precomposition with $\varphi_p$. A direct computation on the Tate curve, whose formal group is canonically $\widehat{\bbG}_m$, shows the cover is unramified at the cusps and thus extends canonically to a finite \'{e}tale $(\bbZ/p\bbZ)^\times$-cover $\Ig^{\ord,c}_{K^p,\mu_p}$ of $\Xint^{\ord}_{K^p, \bbF_p}$.

\begin{definition}[Mod $p$ modular forms]\label{def.mod-p-modular-forms}
The space of mod $p$ modular forms is the $(\bbZ/p\bbZ)^\times \times \GL_2(\bbA_f^{(p)})$-equivariant ring $\mc{M}_{\bbF_p}:=H^0(\Ig^{\ord, c}_{\mu_p}, \mc{O}).$
\end{definition}

The pullback $(\varphi_p^{-1})^*\frac{dt}{t}$ of the invariant differential $\frac{dt}{t}$ on $\mu_p=\Spec \bbF_p[t]/(t^p-1)$ is $\GL_2(\bbA_f^{(p)})$-equivariant trivialization of $\omega$ on $\Ig^{\ord}_{\mu_p}$ and extends to a trivialization over $\Ig^{\ord, c}_{\mu_p}$. In particular, we can use it to evaluate modular forms to functions on $\Ig^{\ord, c}_{\mu_p}$ and we obtain the following version of a well-known result (cf., e.g., \cite{gross:tameness}):

\begin{lemma}\label{lemma.mod-p-modular-forms} Evaluation along $(\varphi_p^{-1})^*\frac{dt}{t}$ induces a $(\bbZ/p\bbZ)^\times \times \GL_2(\bbA_f^{(p)})$-equivariant isomorphism of rings
\[ \left( \bigoplus_{k \geq 0} M_{k,\bbF_p} \right) / (\Ha - 1)	 \cong \mc{M}_{\bbF_p}. \]
where $(\bbZ/p\bbZ)^\times$ acts on  $M_{k,\bbF_p}$ by the character $z \mapsto z^{k}.$ 
\end{lemma}
\begin{proof} 
First note that the evaluation map on $M_{k,\bbF_p}$ factors through $H^0(\Xint^\ord_{\bbF_p},\omega^k)$. One computes that $((\varphi_p^{-1})^*\frac{dt}{t})^{p-1}=\Ha$, so that the evaluation map factors as
\[ \bigoplus_{k\geq 0} M_{k,\bbF_p} \rightarrow \bigoplus_{k=0}^{p-2} H^0(\Xint^\ord_{\bbF_p},\omega^k) \xrightarrow{\sim} H^0(\Ig^{\ord, c}_{\\\mu_p}, \mc{O}) = \mc{M}_{\bbF_p} \]
where the first arrow is restriction followed by division by $\Ha^{\lfloor k/(p-1)\rfloor}$. The second arrow is an isomorphism because we can decompose $H^0(\Ig^{\ord, c}_{\mu_p}, \mc{O})$ according to the characters of $(\bbZ/p\bbZ)^\times$ in $\bbF_p$ and if $f$ is in the character space for $z^{k}$ then $f( (\varphi_p^{-1})^*\frac{dt}{t})^k$ is invariant under $(\bbZ/p\bbZ)^\times$ thus descends to a section over $\Xint^\ord_{\bbF_p}$.

The first map is clearly surjective since multiplying by a sufficient power of $\Ha$ will clear any poles along the supersingular locus. On the other hand, any element in the kernel can be multiplied by $1+\Ha + \Ha^2 + \Ha^3 + \ldots$ (the condition of being in the kernel, when written out, guarantees that this product is zero in sufficiently large degree) to show that it is in fact in the ideal generated by $(1-\Ha)$. 
\end{proof}

\begin{remark} The ideal $(\Ha-1)$ on the left can also be interpreted as the kernel of the total $q$-expansion map, and indeed the original definition of mod $p$ modular forms was by passing to $q$-expansions.  \end{remark}

We now define $p$-adic modular forms. 
\begin{definition}[The Katz-Igusa moduli problem] Let $K^p \leq \GL_2(\bbA_f^{(p)})$ be a closed subgroup. Let $\Ig^{\ord}_{K^p, \wh{\bbG}_m}$ be the functor on $\Nilp_{\bbZ_p}$ 
\[ \Ig^{\ord}_{K^p, \wh{\bbG}_m}: R \mapsto \{(E, \varphi_p, \mc{K}^p)\} / \sim\]
sending $R/\bbZ_p$ to the equivalence classes of triples $(E, \varphi_p, \mc{K}^p)$ such that
\begin{enumerate}
\item $E/R$ is an elliptic curve,
\item $\varphi_p: \wh{\bbG}_m \xrightarrow{\sim} \wh{E}$ is an isomorphism of formal groups
\item $\mc{K}^p \subset \ul{\Isom}( (\ul{\bbA_f^{(p)}})^2, V_{\bbA_f^{(p)}}(E) )$ is a $K^p$-torsor, and
\item The relation $\sim$ is defined by $(E, \varphi_p, \mc{K}^p) \sim (E', \varphi_p', {\mc{K}^p}')$ if there is a prime-to-$p$ quasi-isogeny $q: E \rightarrow E'$ such that $q\circ \varphi_p = \varphi_p'$ and $q(\mc{K}^p)=\mc{K}^p$.	
\end{enumerate}
\end{definition}

\textbf{As usual, when $K^p=\{e\}$ we  drop it from the notation.} The main source is \cite{katz:higher-congruences}, which works with the geometrically connected variant with full level $n$ structure. As in loc. cit., for $K^p$ sufficiently small, $\Ig^{\ord}_{K^p, \wh{\bbG}_m}$ is represented by an affine $p$-adic formal scheme, a pro-\'{e}tale $\bbZ_p^\times$-torsor over the formal ordinary locus $\Yint_{K^p}^{\wedge, \ord}/\Spf\bbZ_p$ -- here the action of $\ul{\bbZ_p^\times}$ is by precomposition with $\varphi_p$. For the same reason as the $\mu_p$-Igusa moduli problem, it is unramified at the cusps and thus extends canonically over $\Xint^{\wedge,\ord}_{K^p} /\Spf \bbZ_p$ to an affine formal scheme $\Ig^{\ord, c}_{K^p, \wh{\bbG}_m}$. 

\begin{definition} The space of $p$-adic modular forms is the unitary $\bbZ_p^\times \times \GL_2(\bbA_f^{(p)})$ representation on the $\bbQ_p$-Banach space
\[ \bbV_{\bbQ_p} := \bbV_{\bbZ_p}[1/p] \textrm{ for } \bbV_{\bbZ_p}:=H^0(\Ig^{\ord,c}_{\wh{\bbG}_m}, \mc{O}). \]
In the introduction we used the notation $\mc{M}_{p-adic}$, which we will not employ further. We note also that $\bbV_{\bbZ_p}$ is $p$-torsion free. 
\end{definition}

The $\GL_2(\bbA_f^{(p)})$-representation $\bbV_{\bbQ_p}$ is not smooth, but one can check that the smooth vectors are dense and, for $K^p$ a sufficiently small closed subgroup, already at the integral level we have
\[ \bbV^{K^p}_{\bbZ_p} = H^0(\Ig^{\ord, c}_{K^p, \wh{\bbG}_m}, \mc{O}). \]

The bundle $\omega$ is $\GL_2(\bbA_f^{(p)})$-equivariantly trivialized over $\Ig^{\ord}_{\wh{\bbG}_m}$ by $(\varphi_p^{-1})^*\frac{dt}{t}$ and this trivialization extends to $\Ig^{\ord,c}_{\wh{\bbG}_m}$. We can thus evaluate modular forms to elements of $\bbV_{\bbQ_p}$ and, as in \cite[especially Theorem 2.1]{katz:higher-congruences}, one obtains

\begin{lemma}\label{lemma:density-p-adic-mf} Evaluation on $(\varphi_p^{-1})^*\frac{dt}{t}$ induces a $\bbZ_p^\times \times \GL_2(\bbA_f^{(p)})$-equivariant injection 
\[ \bigoplus_{k \geq 0} M_{k,\bbQ_p}^{\GL_2(\bbZ_p)} \hookrightarrow \bbV_{\bbQ_p} \]
where we let $\bbZ_p^\times$ act by $z^{k}$ on $M_{k, \bbQ_p}$. The induced map on $K^p$-invariants has dense image for any compact open subgroup $K^p$. 
\end{lemma}

\begin{remark}\label{remark.mod-p-and-p-adic-mf} If we write $\bbV_{\bbF_p}=\bbV_{\bbZ_p}/(p)$, then it is clear by comparing moduli problems that the invariants $\bbV_{\bbF_p}^{1+p\bbZ_p}$ represent $\Ig^{\ord,c}_{\mu_p}.$ The obvious diagram comparing Lemmas \ref{lemma.mod-p-modular-forms} and \ref{lemma:density-p-adic-mf} then commutes: 
\[ \xymatrix{ \bigoplus_{k\geq 0} M_{k, \bbQ_p}^{\GL_2(\bbZ_p)}\ar[d] & \ar[l] \bigoplus_{k\geq 0} H^0(\Xint_{\bbZ_p}, \omega^k)\ar[d]\ar[rr] & & \bigoplus_{k \geq 0} M_{k, \bbF_p}\ar[d] \\ 
\bbV_{\bbQ_p} & \ar[l] \bbV_{\bbZ_p}\ar[r] &  \bbV_{\bbF_p} & \ar[l] \bbV_{\bbF_p}^{1+p\bbZ_p}  = H^0(\Ig^{\ord, c}_{\mu_p}, \mc{O})}\]
This diagram summarizes why in \cite{katz:higher-congruences} one must pass to \emph{divided} congruences to see all of $\bbV$, which in our presentation corresponds to the fact that it is crucial to invert $p$ to obtain the density statement in Lemma \ref{lemma:density-p-adic-mf}. 
\end{remark}

\begin{remark} For $K^p$ a sufficiently small compact open, if we fix a set of cusps $c_1, \ldots, c_m$, one in each connected component of $X_{\GL_2(\ZZ_p)K^p, \breve{\QQ}_p}$, then we obtain $q$-expansion maps
\[ M_{k, \QQ_p}^{\GL_2(\ZZ_p)K^p} \rightarrow \prod_{i=1}^{m} \breve{\ZZ}_p[[q]][1/p]. \]
A generalization of Serre's \cite{serre:formes-modulaires-et-zeta} original definition of $p$-adic modular forms of level $K^p$ would be to take the completion of the span of the images of these maps over all $k$. The $q$-expansion principle combined with Lemma \ref{lemma:density-p-adic-mf} implies that this agrees with the definition given above (see \cite{katz:higher-congruences} for related discussions). 
\end{remark}

\subsection{The supersingular Caraiani-Scholze Igusa variety}\label{ss.igusa-variety}

We now turn our attention to the main player in our story. Though it would be possible to work over $\bbF_p$, from here on out it will be cleaner and more convenient to work over $\overline{\bbF}_p$ (essentially because the endomorphisms of $\bbX_{\ssing, \overline{\bbF}_p}$ are  not all defined over $\bbF_p$).  
 
\begin{definition}[The supersingular Caraiani-Scholze Igusa moduli problem]\label{def.igusa} Let $K^p \leq \GL_2(\bbA_f^{(p)})$ be a closed subgroup. Let $\Ig_{K^p}^\ssing$ be the functor on $\Fpbar$-algebras
\[ \Ig_{K^p}^{\ssing}: R \mapsto  \{(E, \varphi_p, \mc{K}^p)\}/\sim \]
sending $R/\Fpbar$ to the set of equivalence classes of triples 
$(E, \varphi_p, \mc{K}^p )$ where
\begin{enumerate}
\item $E/R$ is an elliptic curve,
\item $\varphi_p: \bbX_{\ssing, R} \otimes \bbQ_p \xrightarrow{\sim} E[p^\infty] \otimes \bbQ_p $ is a quasi-isogeny,
\item $\mc{K}^p \subset \ul{\Isom}\left((\ul{\bbA_f^{(p)}})^2, V_{\bbA_f^{(p)}}(E)\right)$ is a $K^p$-torsor, and
\item the relation $\sim$ is defined by $(E,\varphi_p, \mc{K}^p) \sim (E', \varphi_p', {\mc{K}^p}')$ if there is a quasi-isogeny $q: E \otimes \bbQ \xrightarrow{\sim} E'\otimes \bbQ$ such that 
\[ q(\mc{K}^p) = {\mc{K}^p}' \textrm{ and } q \circ \varphi_p = \varphi'_p. \]
\end{enumerate}
\end{definition}

\textbf{As usual, when $K^p=\{e\}$ we drop it from the notation and write $\Ig^\ssing := \Ig_{\{e\}}^\ssing$}. In this case we  will also identify the prime-to-$p$ level data with the choice of an isomorphism $\varphi_{\bbA_f^{(p)}}: \ul{(\bbA_f^{(p)})^2} \rightarrow V_{\bbA_f^{(p)}}(E)$ as in Example \ref{example.integral-infinite-level}. 

This up to quasi-isogeny moduli problem is well-adapted for comparison to the quaternionic coset, as we will see in \S\ref{ss.quat-coset} below. To compare with modular curves, however, it is useful to also have an up to prime-to-$p$ quasi-isogeny moduli interpretation. Similarly to Lemma \ref{lemma:generic-integral-fiber} (see also \cite[Lemma 4.3.4]{caraiani-scholze:generic}), we obtain:
\begin{lemma}\label{lem.prime-to-p-igusa} For $K^p \leq \GL_2(\bbA_f^{(p)})$ a closed subgroup, consider the functor 
\[  \Ig^\ssing_{K^p, \bbX_\ssing}: R \mapsto  \{(E, \varphi_p, \mc{K}^p)\}/\sim \]
sending $R/\Fpbar$ to the set of equivalence classes of triples 
$(E, \varphi_p, \mc{K}^p)$ where
\begin{enumerate}
\item $E/R$ is an elliptic curve,
\item $\varphi_p:  \bbX_{\ssing,R} \xrightarrow{\sim} E[p^\infty] $  is an \textbf{isomorphism,}
\item $\mc{K}^p \subset \ul{\Isom}\left((\ul{\bbA_f^{(p)}})^2, V_{\bbA_f^{(p)}}(E)\right)$ is a $K^p$-torsor, and
\item the relation $\sim$ is defined by $(E,\varphi_p, \mc{K}^p) \sim (E', \varphi_p', {\mc{K}^p}')$ if there is a \textbf{prime-to-$p$ quasi-isogeny} $q: E \otimes \bbZ_{(p)} \xrightarrow{\sim} E'\otimes \bbZ_{(p)}$ such that 
\[ q(\mc{K}^p) = {\mc{K}^p}' \textrm{ and } q \circ \varphi_p = \varphi'_p. \]
\end{enumerate}
The assignment 
\[ (E, \varphi_p, \mc{K}^p) \mapsto (E, \varphi_p \otimes \bbQ_p, \mc{K}^p) \]
induces an isomorphism 
\[ \Ig^\ssing_{K^p, \bbX_\ssing} \rightarrow \Ig^\ssing_{K^p}. \]
\end{lemma}

We write $\mc{O}_p = \End(\bbX_{\ssing, \Fpbar})$ and $D_p := \mc{O}_p \otimes \mbb{Q}_p$, so that $\mc{O}_p$ is the maximal order in $D_p$, a ramified quaternion algebra over $\bbQ_p$. Then $D_p$ acts on $\mbb{X}_{\ssing, \Fpbar} \otimes \bbQ_p$, and there is a natural action of $D_p^\times$ on $\Ig^\ssing_K $ by composition with $\varphi_p$. The action of $\mc{O}_p^\times \subset D_p^\times$ preserves the prime-to-$p$ moduli interpretation of Lemma \ref{lem.prime-to-p-igusa}. 

In fact, what acts most naturally are the functors on $\Fpbar$-algebras
\[ \ul{\Aut}(\mbb{X}_{\ssing, \Fpbar}): R \mapsto \Aut(\mbb{X}_{\ssing, R}) \textrm{ and } \ul{\Aut}(\mbb{X}_{\ssing, \Fpbar} \otimes \bbQ_p): R \mapsto \Aut (\mbb{X}_{\ssing, R} \otimes \bbQ_p). \]
The following lemma says we have not missed anything:
\begin{lemma}\label{lemma.iso-isog-p-div} The actions of $D_p^\times$ and $\mc{O}_p^\times$ described above factor through isomorphisms 
	$\ul{\mc{O}_p^\times} \rightarrow \ul{\Aut}(\mbb{X}_{\ssing,\Fpbar})$ and $\ul{D_p^\times} \rightarrow \ul{\Aut}(\mbb{X}_{\ssing,\Fpbar} \otimes \bbQ_p)$. 
\end{lemma}
\begin{proof} This is established in a more general context in the proof of \cite[Proposition 4.2.11]{caraiani-scholze:generic}. The key point is that on $\overline{\FF}_p$-algebras, by \cite[Lemma 4.1.7, Corollary 4.1.10]{caraiani-scholze:generic}, the endomorphisms of $\mbb{X}_{\ssing, \overline{\bbF}_p}$ are given by the Tate module of an \'{e}tale $p$-divisible group over ${\overline{\bbF}_p}$, and are thus equal to the constant sheaf on the ${\overline{\bbF}_p}$-points of that Tate module.  
\end{proof}

\begin{remark} One can also define an ordinary Caraiani-Scholze Igusa variety. We will not use this here, but see Remark \ref{remark.final-caraiani-scholze} for some connections to the present work and \cite{howe:unipotent-circle} for an application to the study of $p$-adic modular forms. \end{remark}

\subsection{Uniformization of the supersingular locus}\label{ss.uniformization}

\begin{proposition}\label{prop.mod-p-unif}
For $K^p$ a sufficiently small closed subgroup, the assignment 
\[ (E, \varphi_p, \mc{K}^p ) \mapsto (E, \mc{K}^p ) \]
induces a $\ul{D_p^\times \times N_{\GL_2(\mbb{A}_f^{(p)})}(K^p)}$-equivariant map
\[ \unif_{\bbX_\ssing}: \Ig^\ssing_{K^p}= \Ig^{\ssing}_{K^p,\bbX_\ssing} \rightarrow \Yint^\ssing_{K^p, \overline{\bbF}_p}, \]
where $\ul{D_p^\times}$ acts on $\Yint^\ssing_{K^p, \overline{\bbF}_p}$ by $\Frob^{v_p\circ \Nrd}$ (here $\Frob: \Yint_{K^p,\Fpbar} \rightarrow \Yint_{K^p,\Fpbar}$ is the relative Frobenius and $\Nrd: D_p \rightarrow \bbQ_p$ denotes the reduced norm). 
Moreover, $\unif_{\bbX_{\ssing}}$ is a trivializeable $\ul{\calO_p^\times}$-torsor, and the maps compile to a $\GL_2(\bbA_f^{(p)})$-equivariant isomorphism of towers as the sufficiently small closed subgroup $K^p$ varies.  
\end{proposition}
\begin{proof}
We first show that the map factors through $\Yint_{K^p,\overline{\bbF}_p}^\ssing$: given a point 
\[ (E, \varphi_p, \mc{K}^p) \in \Ig^\ssing_{K^p,\bbX_\ssing}(R),\]
we find $\Ha(E/R)$ is identically zero because $E[p^\infty] \cong \bbX_{\ssing,R}$ and, as observed in \S\ref{sss.hasse}, the Hasse invariant only depends on $E[p^\infty]$. But $\Ha(E/R)$ is the pullback of $\Ha$ under the induced map $\Spec R \rightarrow \Yint_{K^p, \Fpbar}$, so by Lemma \ref{lemma.Hasse-vanishing} this induced map factors through $\Yint^\ssing_{K^p, \Fpbar}$, as desired.  

It is also clear from the definitions of the moduli problems and the  identification $\ul{\Aut}(\bbX_{\ssing, \Fpbar})=\ul{\calO_p^\times}$ in Lemma \ref{lemma.iso-isog-p-div} that the map is a quasi-torsor; to conclude it is a trivializeable torsor it thus suffices to produce a section. 

To obtain this, we observe that the $p$-divisible group of the universal elliptic curve over $\Yint_{K^p, \Fpbar}^\ssing$ is isomorphic to $\bbX_{\ssing, \Yint_{K^p, \Fpbar}^\ssing}$: indeed, it is pulled back from $\Yint_{{K^p}', \Fpbar}^\ssing$ for any sufficiently small compact open subgroup ${K^p}' \leq \GL_2(\bbA_f^{(p)})$ containing $K^p$,  and $\Yint_{{K^p}', \Fpbar}^\ssing$ is just a finite union of $\overline{\bbF}_p$-points. We can choose an isomorphism at each of these points (because $\bbX_{\ssing, \overline{\bbF}_p}$ is the unique up to isomorphism connected one dimensional height two $p$-divisible group over $\overline{\bbF}_p$), and then assemble these and pull back to $\Yint_{K^p}^\ssing$ to get the desired isomorphism. This is exactly the data of a section. 

It is left only to verify the action of $\ul{D_p^\times}/\ul{\mc{O}_p^\times}$ is as described on $\Yint_{K^p,\Fpbar}^\ssing$, and this is a direct computation from the definitions that we leave to the reader (we will not use this part of the statement in any of what follows). 
\end{proof}

As a consequence, one finds that for each sufficiently small $K^p$, $\Ig_{K^p}^\ssing$ is a profinite set, thus in particular a perfect affine scheme over $\overline{\mbb{F}}_p$. Our next goal is to identify this profinite set explicitly with a quaternionic coset.

\subsection{Supersingular Igusa variety as quaternionic coset}\label{ss.quat-coset}
We will now identify $\Ig^\ssing$ with a quaternionic coset using the orbit map for the action of $\ul{D_p^\times \times \GL_2(\mbb{A}_f^{(p)})}$. To that end, fix a supersingular elliptic curve $E_0 / \overline{\mbb{F}}_p$ and level structure to obtain 
\[ x_0 = (E_0, \varphi_{p,0}, \varphi_{\bbA_f^{(p)},0}) \in \Ig^\ssing({\overline{\bbF}_p}). \] 

We write $D=\End(E_0) \otimes \mbb{Q}$. As explained, e.g., in \cite[Theorem V.3.1]{silverman:aoec}, this a quaternion algebra over $\bbQ$, and because it is nonsplit and acts faithfully on the $\ell$-adic Tate module for all $\ell \neq p$, it must be ramified exactly at $p$ and $\infty$. 

By definition, $D^\times$ is identified with the self quasi-isogenies of $E_0 \otimes \mbb{Q}$. The actions of $D$ on $E_0[p^\infty]\otimes\bbQ_p$ and $V_{\bbA_f^{(p)}}(E)$, transported by $\varphi_{p,0}$ and $\varphi_{\bbA_f^{(p)},0}$, thus induce an identification 
\[ D^\times(\bbA_f) = D^\times(\bbQ_p) \times D^\times(\bbA_f^{(p)})= D_p^\times \times \GL_2(\bbA_f^{(p)}). \]
In particular, we obtain a map $D^\times(\bbQ) \hookrightarrow D^\times(\bbA_f)$. 

The action on the point $x_0$ induces an orbit map 
\begin{equation}\label{eq.orbit-map} \underline{D_p^\times \times \GL_2(\bbA_f^{(p)})} \rightarrow \Ig^\ssing, \; \;  g \mapsto x_0 g, \end{equation} 
and we will show: 
\begin{theorem}\label{theorem:orbit-coset} The orbit map (\ref{eq.orbit-map}) factors through a ${D^\times(\bbA_f)=D^\times(\bbQ_p) \times \GL_2(\bbA_f^{(p)})}$-equivariant isomorphism 
\[ \underline{ D^\times(\bbQ) \bs D^\times(\bbA_f)} = \underline{D^\times(\bbQ) \bs \l( D_p^\times \times \GL_2\l(\bbA_f^{(p)}\r) \r)} \xrightarrow{\sim} \Ig^{\ssing}. \]
\end{theorem}

Before proving Theorem \ref{theorem:orbit-coset}, it will be helpful to recall the following basic structural result for quaternionic double cosets. We write $\calO = \End(E_0)$, a maximal order in $D$, and consider the \emph{finite} class set of right fractional ideals in $\calO$ 
\[ D^\times(\bbQ) \bs D^\times(\bbA_f) / (\calO \otimes \wh{\bbZ})^\times = D^\times(\bbQ) \bs D_p^\times \times \GL_2(\bbA_f^{(p)}) / \calO_p^\times \times \GL_2(\wh{\bbZ}^{(p)}). \]
We may fix a finite set of representatives for the class set $I$ and corresponding representatives $\gamma_\calI \in D^\times(\bbA_f),\, \calI \in I$ for the double cosets. The stabilizer of $\calI$ for left multiplication in $D^\times(\bbQ)$ is the units in a maximal order $\calO_\calI$, and we find 
\[ D^\times(\bbQ) \bs D^\times(\bbA_f)  \cong  \bigsqcup_{\calI \in I} \calO_\calI^\times \bs \gamma_\calI (\calO \otimes \wh{\bbZ})^\times = \bigsqcup_{\calI \in I} \calO_{\calI}^\times \bs \Isom(\calO \otimes \wh{\bbZ}, \calI \otimes \wh{\bbZ}) \]
where the the isomorphisms on the right are of right $\calO \otimes \wh{\bbZ}$-modules. Because each group $\calO_\calI^\times$ is finite, if we replace $\calO \otimes \wh{\bbZ}$ with a small enough compact open subgroup $K \subset \calO \otimes \wh{\bbZ}^\times$, we obtain a finite set of representatives $\gamma_{\calI, i}$ such that
\begin{equation} D^\times(\bbQ) \bs D^\times(\bbA_f) \cong \bigsqcup_{\calI \in I, i} \gamma_{\calI,i} K \end{equation}
In other words, we obtain a topological splitting of the locally profinite set $D^\times(\bbA_f)$ as a product of a discrete set and a profinite set,
\begin{equation}\label{equation:decomp-quat}  D^\times(\bbA_f) \cong D^\times(\bbQ) \times D^\times(\bbQ) \bs D^\times(\bbA_f) \end{equation}
compatible with the left action of $D^\times(\bbQ)$ and the right action of $K$. 

We will also need to understand the quasi-isogenies of $E_0$ after arbitrary base change. To this end, we observe that, for any $S$, there is a natural map from $\ul{D^\times(\bbQ)}(R)$ to  $\ul{\Aut}(E_0\otimes \bbQ)(R)=\Aut(E_{0,R} \otimes \bbQ)$ -- indeed, since $D^\times(\bbQ)$ is discrete, an element of $\ul{D^\times(\bbQ)}(R)$ is a locally constant on $\Spec R$ choice of element in $D^\times(\bbQ)$. 
\begin{lemma}\label{lemma.self-qisog-simple} The natural map defined above identifies $\ul{D^\times(\bbQ)}$ with $\ul{\Aut}(E\otimes \bbQ)$
\end{lemma}
\begin{proof}
	It suffices to verify this on $\Spec R$ points, and we may moreover assume $R$ is of finite type over ${\overline{\bbF}_p}$, and then that $R$ is reduced using that quasi-isogenies lift along nilpotent ideals containing $p$. Then, the result follows from the computation over algebraically closed fields. 
\end{proof}

\begin{proof}[Proof of Theorem \ref{theorem:orbit-coset}]
By Lemma \ref{lemma.self-qisog-simple}, we deduce that the orbit map factors as an injection on $R$-points
\[ \underline{D^\times(\bbQ)}(R) \bs \underline{D^\times(\bbA_f)}(R) \rightarrow \Ig^\ssing (R). \]
From (\ref{equation:decomp-quat}) we deduce 
\[ \underline{D^\times(\bbQ) \bs D^\times(\bbA_f)} (R) = \underline{D^\times(\bbQ)}(R) \bs \underline{D^\times(\bbA_f)}(R), \]
and thus the orbit map factors through an injection 
\[ \underline{D^\times(\bbQ) \bs D^\times(\bbA_f)} \hookrightarrow 
\Ig^{\ssing}. \]

It remains to show the map is surjective. To do so, it suffices to show that the universal elliptic curve over $\Ig_\ssing$ is quasi-isogenous to $E_{0,\Ig^\ssing}$. But the universal elliptic curve is the pullback from $Y^\ssing_{\Fpbar}$ along the map $\unif_{\bbX_\ssing}$ of Proposition \ref{prop.mod-p-unif}, and the statement follows as in the proof of Proposition \ref{prop.mod-p-unif} by reduction to a finite set of points at finite level and the fact that any two supersingular elliptic curves over ${\overline{\bbF}_p}$ are isogenous.  
\end{proof}

\begin{remark} We can avoid the use of any Grothendieck topology above because our torsors are all trivializable. In particular, this sidesteps the following question: for $X$ a topological space and $\Gamma$ a topological group acting on $X$, for which Grothendieck topologies does $\ul{X}/\ul{\Gamma}=\ul{X/\Gamma}$? Indeed, when $X \rightarrow X/\Gamma$ is a trivializable $\Gamma$-torsor, this is already true at the level of presheaves.\end{remark}

\begin{corollary}\label{corollary.igusa-level-coset} For any sufficiently small closed subgroup $K^p\leq \GL_2(\bbA_f^{(p)})$, the map of Theorem \ref{theorem:orbit-coset} induces a $\ul{D_p^\times \times N_{\GL_2(\bbA_f^{(p)})}(K^p)}$-equivariant isomorphism  
\[ \underline{D^\times(\bbQ) \bs \l( D_p^\times \times \GL_2\l(\bbA_f^{(p)}\r) \r)/K^p }  \xrightarrow{\sim}  \Ig^{\ssing}_{K^p}. \]
and, combined with Proposition \ref{prop.mod-p-unif}, 
\[ \ul{D^\times(\bbQ) \bs D^\times(\bbA_f) / \l(\calO_p^\times \times K^p\r)}\xrightarrow{\sim} \Yint_{K^p, \overline{\bbF}_p}^\ssing. \]
These maps compile to $\GL_2(\bbA_f^{(p)})$-equivariant isomorphisms of towers as the closed subgroup $K^p$ varies over all sufficiently small closed subgroups. 
\end{corollary}

\section{Serre's mod $p$ correspondence}\label{section:mod-p}

In this section we give a proof of Theorem \ref{theorem:serre}, roughly following Serre \cite{serre:two-letters}. The main difference between our presentation and that of loc. cit. is that we emphasize from the beginning the role of the supersingular Igusa variety as quaternionic coset in the uniformization of the supersingular locus established in Corollary \ref{corollary.igusa-level-coset}. When $p=2 \textrm{ or } 3$ the proof is only valid for $K^p$ sufficiently small, however, in \S\ref{ss.deduction-mod-p-from-p-adic} we explain how to deduce the full result directly from Theorem \ref{mainthm:HeckeIso}. 

In this section $D/\bbQ$ is the specific quaternion algebra ramified at $p$ and $\infty$ defined in \S\ref{ss.quat-coset}. Up to isomorphism, there is a unique quaternion algebra over $\bbQ$ ramified only at $p$ and $\infty$, so this specific choice is made only to normalize with the results of \S\ref{ss.quat-coset}. Recall that we have also fixed in \S\ref{ss.quat-coset} identifications $D^\times(\bbA_f^{(p)})=\GL_2(\bbA_f^{(p)})$ and $D^\times(\bbQ_p)=D_p^\times$. As in the introduction we write
\[ \mc{A}_{\Fpbar} = \Cont(D^\times(\bbQ) \bs D^\times(\bbA_f), \Fpbar), \]
which we view as an admissible  $D^\times(\bbA_f)=D_p^\times \times \GL_2(\bbA_f^{(p)})$-representation.
\subsection{Mod $p$ modular forms}

We write $\mc{M}_{\Fpbar}$ for the base change to $\overline{\bbF}_p$ of the space of mod $p$ modular forms of Definition \ref{def.mod-p-modular-forms}. As explained in the proof of Lemma \ref{lemma.mod-p-modular-forms}, the evaluation map from modular forms to mod $p$ modular forms factors through an isomorphism of $\GL_2(\bbA_f^{(p)})$-representations 
\begin{equation}\label{eq.sections-ordinary-locus} \mc{M}_{\Fpbar} = \bigoplus_{k=0}^{p-2} H^0(\frak{X}_{\Fpbar}^\ord, \omega^k). \end{equation}
We consider the increasing exhaustive filtration $F_i \mc{M}_{\overline{\bbF}_p}$ whose $i$th step consists of those sections on the right of (\ref{eq.sections-ordinary-locus}) with poles of order $\leq i$ along $\mathfrak{X}_{\Fpbar}^\ssing$.  

As in \S\ref{ss.modular-forms} we also write 
\[ M_{k, \Fpbar} = H^0(\Xint_{\Fpbar}, \omega^k), \]
the admissible $\GL_2(\bbA_f^{(p)})$ representation of weight $k$ modular forms over $\overline{\bbF}_p$. In particular, for $K^p$ sufficiently small, we have
\[ M_{k, \Fpbar}^{K^p} = H^0(\Xint_{K^p, \Fpbar}, \omega^k). \]

From these definitions and Lemma \ref{lemma.Hasse-vanishing}, we find that multiplication by $\Ha^i$ gives a $\GL_2(\bbA_f^{(p)})$-equivariant isomorphism 
\begin{equation}\label{eq.mod-p-filt-iso} F_i \mc{M}_{\overline{\bbF}_p} \xrightarrow{\sim} \bigoplus_{k=0}^{p-2} M_{k+(p-1)i,\Fpbar}. \end{equation} 
Thus, in particular, $F_i \mc{M}_{\overline{\bbF}_p}$ is an admissible $\GL_2(\bbA_f^{(p)})$-representation.

\subsection{Evaluation of modular forms} We now explain how to evaluate modular forms to functions on $\Ig^\ssing$. This works almost exactly as in the evaluation map for mod $p$ modular forms in Lemma \ref{lemma.mod-p-modular-forms}, but we will set things up here in a slightly more canonical language.

We write $\omega_{\bbX_\ssing}=\Lie(\bbX_\ssing)^*$, equipped with the natural action of $\calO_p^\times$ after base change to $\Fpbar$. Using the isomorphism $\varphi_p$ in the prime-to-$p$ moduli interpretation (cf. Definition \ref{def.igusa}-(2)) and the uniformization in Proposition \ref{prop.mod-p-unif}, 
 we obtain a canonical $\calO_p^\times \times \GL_2(\bbA_f^{(p)})$-equivariant isomorphism 
 \[ \omega_{\bbX_\ssing} \otimes_{\bbF_p} \calO_{\Ig^{\ssing}} \xrightarrow{\sim} \unif_{\bbX_\ssing}^* \omega. \]
   In particular, we obtain a $\GL_2(\bbA_f^{(p)})$-equivariant isomorphism
\[ H^0(\Xint_{\Fpbar}^\ssing, \omega^k) \xrightarrow{\sim} \Hom_{\calO_p^\times}( \Lie(\bbX_{\ssing,\Fpbar})^k, H^0(\Ig^\ssing, \calO) ). \]

Note that, by the definition of $\bbX_\ssing$ in terms of the Dieudonn\'{e} module given in \ref{ss.igusa-variety}, we have an identification of $\Lie \bbX_{\ssing}$ with $\bbF_p^2/\langle (1, 0) \rangle$, and of $\calO_p$ with the $\sigma$-centralizer of $b_\ssing$ in $M_2(\breve{\bbZ}_p)$. In particular, the image of $(0,1)$ gives a basis element of $\Lie \bbX_\ssing$, and a direct computation shows $\calO_p^\times$ acts on $\Lie \bbX_{\ssing,\Fpbar}={\overline{\bbF}_p}$ through a surjective character $\epsilon: \mc{O}_p^\times \rightarrow \mbb{F}_{p^2}^\times$ whose kernel is a pro-$p$ group we write as $\calN_p$. Combining with Theorem \ref{theorem:orbit-coset}, we obtain a $\GL_2(\bbA_f^{(p)})$-equivariant isomorphism
\begin{equation}\label{eqn:mod-p-eq-ev} H^0(\Xint_{\Fpbar}^\ssing, \omega^k) \xrightarrow{\sim} \Hom_{\calO_p^\times}\l( \epsilon^k, \Cont\l(D^\times(\bbQ) \bs D_p \times \GL_2(\bbA_f^{(p)}) / \calN_p, {\overline{\bbF}_p} \r) \r). \end{equation}

Then, evaluating homomorphisms in the right-hand side of (\ref{eqn:mod-p-eq-ev}) on $1 \in {\overline{\bbF}_p}$ and passing to $K^p$-invariants, we obtain Serre's map, a Hecke equivariant isomorphism
\begin{equation}\label{equation:serre-map}  H^0(\Xint_{\Fpbar}^\ssing, \omega^k) \xrightarrow{\sim} \calA^{\calN_p K^p}_{{\overline{\bbF}_p}}[\epsilon^k]. \end{equation}

\subsection{Hecke algebras and generalized eigenspaces}\label{ss.generalized-eigenspaces}
If $L$ is a field and $V$ is a finite dimensional vector space over $L$, then for any commutative $L$-algebra $T$ acting on $V$, we have a decomposition into generalized eigenspaces 
\[ V = \bigoplus V_{\frakm} \]
where $\frakm$ runs over maximal ideals of $T$ with residue field a finite extension of $L$ and $V_\frakm$ is the sub-module of $\frakm$-torsion elements. In particular, if $V_\frakm \neq 0$, the eigenspace $V[\frakm]$ consisting of elements annihilated by $\frakm$ is non-empty.

The same result applies more generally to $V$ an increasing union of finite dimensional vector spaces with $T$-action. In particular, if $W$ is an increasing union of admissible representations of a locally profinite group $G$, $K$ is a compact open subgroup of $G$, and $T$ is a commutative subalgebra of the abstract Hecke algebra $L[K\bs G/K]$, then the action of $T$ on the invariants $W^K$ admits a decomposition into generalized eigenspaces. 

This formalism gives a decomposition into generalized eigenspaces for the action of any commutative subalgebra
\[ \bbT' \subset \overline{\bbF}_p[K^p\bs\GL_2(\bbA_f^{(p)}/K^p] = \bbT^\abs_{K^p} \otimes \overline{\bbF}_p, \]
on $M_{k, \Fpbar}^{K^p}$, $\mc{A}^{K^p}_{\Fpbar}$, and $\mc{M}_{\Fpbar}^{K^p}$. For $M_{k,\Fpbar}^{K^p}$ this is immediate because $M_{k, \Fpbar}$ is admissible. For $\mc{A}^{K^p}_{\Fpbar}$ it follows because $\mc{A}$ is a colimit of admissible $\GL_2(\bbA_f^{(p)})$-representations by taking invariants under compact open subgroups of $D_p^\times$. For $\mc{M}_{\Fpbar}^{K^p}$ it follows using the filtration by the admissible representations $F_i \mc{M}_{\Fpbar}.$ 

\begin{remark}\label{remark:mod-p-twist} In order to obtain a $\GL_2(\bbA_f^{(p)})$-action on modular forms inducing the standard action of Hecke operators on $K^p$-invariants, one must twist by the unramified determinant character that appears in the Kodaira-Spencer isomorphism  $\omega^2 = \Omega(\cusps)$. This replaces each individual coset Hecke operators with a multiple by an invertible element (since we do not include Hecke operators at $p$!), so it will not change the image under the action map. Thus this choice is immaterial for our purposes; we prefer not to include the twist because it is less natural to do so except when one is comparing with \'{e}tale cohomology, where it is baked in. 
\end{remark}

\subsection{Spectral decompositions}
We now compare the spectral decompositions provided by the previous section in order to prove Theorem \ref{theorem:serre}. The comparison of $\mc{M}_{\Fpbar}^{K^p}$ and $\mc{A}^{K^p}_{\Fpbar}$ is mediated through comparisons of each with $\bigoplus_{k\geq 0}M_{k, \Fpbar}^{K^p}$.

We first treat the simpler case of $\mc{M}_{\Fpbar}^{K^p}$:

\begin{lemma} \label{lemma:spec-decomp-mf-mod-p}
For $\frakm$ a maximal ideal of $\bbT'$, 
\[ \l( \mc{M}^{K^p}_{\overline{\bbF}_p} \r)_\frakm \neq 0\iff \bigoplus_{k\geq 0}\l(M_{k, \Fpbar}^{K^p}\r)_\frakm \neq 0. \]
\end{lemma} 
\begin{proof}
We have $\l( \mc{M}^{K^p}_{\overline{\bbF}_p} \r)_\frakm \neq 0$ if and only if $\l(F_i \mc{M}^{K^p}_{{\overline{\bbF}_p}} \r)_\frakm \neq 0$ for $i$ sufficiently large. On the other hand, by \ref{eq.mod-p-filt-iso}, $\l(F_i \mc{M}^{K^p}_{{\overline{\bbF}_p}} \r)_\frakm \neq 0$ if and only if $\l(M_{k' + i(p-1)}^{K^p}\r)_\frakm \neq 0$ for some $0 \leq k' \leq p-2$, and we conclude as varying $k'$ and $i$ exhausts all possible $k=k'+i(p-1).$ \end{proof}

On the other hand, using Serre's evaluation map we obtain:
\begin{lemma}\label{lemma:spec-decomp-quat-mod-p} Assume either $K^p$ is sufficiently small or that $p \neq 2, 3$. Then, for $\frakm$ a maximal ideal of $\bbT'$,
\[  \l(\calA^{\calN_p K^p}_{{\overline{\bbF}_p}}\r)_{\frakm} \neq 0 \iff \bigoplus_{k\geq 0} \l(M_{k,\Fpbar}^{K^p}\r)_\frakm \neq 0. \]
\end{lemma}
\begin{proof} 
Suppose $\bigoplus_{k \geq 0} \l(M_{k,\Fpbar}^{K^p}\r)_\frakm \neq 0$. For each $k \geq 0$, we have a $\GL_2(\bbA_f^{(p)})$-equivariant exact sequence
\[ 0 \rightarrow M_{k-(p-1), \Fpbar} \xrightarrow{\cdot H}  M_{k,\Fpbar} \xrightarrow{\mathrm{restriction}} H^0(\Xint_{\Fpbar}^\ssing, \omega^{k}). \]
Passing to $K^p$-invariants and localizing at $\frakm$, we obtain the exact sequence

\begin{equation} \label{equation:mod-p-hasse-localized} 0 \rightarrow \l( M_{k-(p-1), \Fpbar}^{K^p}\r)_\frakm \xrightarrow{\cdot H}  \l(  M_{k,\Fpbar}^{K^p}\r)_\frakm \xrightarrow{\mathrm{restriction}} \l(H^0(\Xint_{\Fpbar}^\ssing, \omega^{k})^{K^p}\r)_\frakm. \end{equation}
By induction on $k$, we deduce that if $\l(M_{k,\Fpbar}^{K^p}\r)_\frakm \neq 0$ then $\left(H^0(\Xint_{\Fpbar}^\ssing, \omega^{k'})^{K^p}\right)_\frakm \neq 0$ for some $0 \leq k' \leq k$, and, applying (\ref{equation:serre-map}), that $\l(\calA^{\mc{N}_pK^p}_{{\overline{\bbF}_p}}\r)_\frakm \neq 0.$ 

Suppose $\l(\calA^{\calN_p K^p}_{{\overline{\bbF}_p}}\r)_\frakm \neq 0.$ Then, as 
\[ \calA^{\calN_pK^p}_{{\overline{\bbF}_p}} = \bigoplus_{k \in \bbZ/(p^2-1)\bbZ}\calA^{\calN_pK^p}_{{\overline{\bbF}_p}}[\epsilon^k], \]
we deduce that for some some $k \geq 0$, $\calA^{\calN_pK^p}_{{\overline{\bbF}_p}}[\epsilon^k]_\frakm \neq 0.$ If $K^p$ is sufficiently small, then because $\omega$ is ample on $\Xint_{K^p,\Fpbar}$ and $\epsilon^k$ only depends on $k$ mod $p^2-1$, we may choose this value of $k$ large enough that ${H^1(\Xint_{K^p,\Fpbar}, \omega^{k-(p-1)})=0}.$ Then, the sequence (\ref{equation:mod-p-hasse-localized}) extends to a short exact sequence; indeed, it is obtained by localizing at $\mathfrak{m}$ the short exact sequence
\[ 0 \rightarrow H^0(\Xint_{K^p,\Fpbar}, \omega^{k-{p-1}}) \xrightarrow{\cdot \Ha} H^0(\Xint_{K^p,\Fpbar}, \omega^{k}) \rightarrow H^0(\Xint_{K^p,\Fpbar}^\ssing, \omega^{k}) \rightarrow  \underbrace{H^1(\Xint_{K^p,\Fpbar}, \omega^{k-(p-1)})}_{=0}.\]
 In particular, restriction induces a surjection
\begin{equation}\label{eq.serre-restriction-surjection} (M_{k,\Fpbar}^{K^p})_\frakm \twoheadrightarrow H^0(\Xint^\ssing_{K^p,\Fpbar}, \omega^{k})_\frakm. \end{equation}
Applying Serre's isomorphism (\ref{equation:serre-map}) we obtain $H^0(\Xint^\ssing_{K^p,\Fpbar}, \omega^{k})_\frakm \neq 0 $, and thus (\ref{eq.serre-restriction-surjection}) implies $(M_{k,\Fpbar}^{K^p})_\frakm \neq 0$, as desired. If $K^p$ is not sufficiently small then we can still apply the same argument by first passing to a sufficiently small $K_1^p$ that is normal in $K^p$ and such that $p$ does not divide $[K^p:K_1^p]$, then taking $K^p/K^p_1$ invariants (which are exact because $p \nmid |K^p/K^p_1|$). When $p\neq 2, 3$, such a $K_1^p$ always exists by adding full level $\ell$ structure for a large enough prime $\ell \not\equiv \pm 1 \mod p$.\end{proof}

Finally a purely representation-theoretic argument allows us to pass back and forth between $\mc{A}_{\Fpbar}^{K^p}$ and its $\mc{N}_p$-invariants: 
\begin{lemma}\label{lemma.d-to-np} For $\frakm$ a maximal ideal of $\bbT'$,
\[ \l(\calA^{K^p}_{{\overline{\bbF}_p}}\r)_{\frakm} \neq 0 \iff \l(\calA^{\calN_p K^p}_{{\overline{\bbF}_p}}\r)_{\frakm} \neq 0 \]
\end{lemma} 
\begin{proof}
Because $\calN_p$ is a pro-$p$ group and $\l(\calA^{K^p}_{{\overline{\bbF}_p}}\r)_\frakm \subset \calA^{K^p}_{{\overline{\bbF}_p}}$ is a smooth characteristic $p$ representation of $\mc{N}_p$, if the representation is nonzero then it admits a nonzero $\calN_p$-fixed vector by the standard trick: the $K_p$-invariants for some $K_p \trianglelefteq \mc{N}_p$ compact open are then non-zero, and then the orbit-stabilizer theorem applied to the action of $\mc{N}_p/K_p$ on the finite dimensional $\bbF_p$-vector space spanned by the orbit of a non-trivial $K_p$-fixed vector shows there must be a non-trivial $\mc{N}_p$-fixed vector.  
\end{proof}

\begin{remark} One can also obtain a similar statement relating $\mc{M}^{K^p}_{\Fpbar}$ and the mod $p$ reduction of $p$-adic modular forms, $\bbV_{\overline{\bbF}_p}^{K^p}$ (see \ref{ss.ordinary-igusa-mod-p-padic} for the notation). Indeed, as in Remark \ref{remark.mod-p-and-p-adic-mf}, $\mc{M}^{K^p}_{\Fpbar}=(\bbV^{K^p}_{\overline{\bbF}_p})^{1+p\bbZ_p}$, and $1+p\bbZ_p$ is a pro-$p$-group; however, a little more work is necessary to obtain the finite dimensionality needed to apply the formalism of the previous section to $\bbV^{K^p}_{\Fpbar}$ in the first place. A version of this argument appears in Lemma \ref{lemma:hecke-action-open-p-adic-mf} below. 
\end{remark}

\subsubsection{Consequences}

Combining Lemmas \ref{lemma:spec-decomp-mf-mod-p}, \ref{lemma:spec-decomp-quat-mod-p}, and \ref{lemma.d-to-np}, we obtain Theorem \ref{theorem:serre} of the introduction with the additional restriction that $K^p$ be sufficiently small when $p=2$ or $3$; only a finite number of maximal ideals can appear in the decompositions because $\calA^{\calN_p K^p}_{\overline{\bbF}_p}$ is a finite dimensional ${\overline{\bbF}_p}$-vector space. We note that this also proves a theorem of Jochnowitz \cite{jochnowitz:congruences} stating that there are only finitely many eigensystems appearing in mod $p$ modular forms. 

\subsection{Deduction of Theorem \ref{theorem:serre} from Theorem \ref{mainthm:HeckeIso}}\label{ss.deduction-mod-p-from-p-adic}
To finish the section, we explain how one can deduce the mod $p$ correspondence (Theorem \ref{theorem:serre}) directly from the $p$-adic correspondence (Theorem \ref{mainthm:HeckeIso}), assuming the latter has been established (as will be done independently of the mod $p$ correspondence in \S\ref{section:p-adic}). In fact, it is more natural to deduce the corresponding statement over $\bbF_p$ instead of $\overline{\bbF}_p$, but there is no important difference between the two.   

We fix a compact open $K^p \leq \GL_2(\bbA_f^{(p)})$ and a commutative subalgebra $\bbT'$ of $\bbT^\abs$. Then, it suffices to verify that the maximal ideals appearing in the decompositions given by the formalism of \S\ref{ss.generalized-eigenspaces} for $\mc{A}_{\bbF_p}^{K^p}$ and $\mc{M}_{\bbF_p}^{K^p}$ are identified with the \emph{open} maximal ideals in the corresponding completed Hecke algebras. We show this in two lemmas, one for each space.

\begin{lemma}\label{lemma:hecke-action-open-dtimes} The action of $\bbT'$ on $\mc{A}_{\bbF_p}^{K^p}$ factors through $\bbT'^\wedge_{\mc{A}_{\bbQ_p}^{K^p}}$ and the maximal ideals such that $\left(\mc{A}_{\bbF_p}^{K^p}\right)_\frakm \neq 0$ are precisely the open ideals of $\bbT'^\wedge_{\mc{A}_{\bbQ_p}^{K^p}}$. 
\end{lemma}
\begin{proof}
The union of the finite dimensional subspaces $\mc{A}_{\bbQ_p}^{K_pK^p}$ as $K_p$ ranges over compact open subgroups $K_p \leq D_p^\times$ is dense in $\mc{A}_{\bbQ_p}^{K^p}$. Thus, if we write $\bbT'_{K_pK^p,n}$ for the image of $\bbT'$ in  $\End_{\bbZ/p^n}(\mc{A}_{\bbZ/p^n\bbZ}^{K_pK^p})$, Lemma \ref{lemma:dense-prodiscrete} gives
\[ \bbT'^\wedge_{\mc{A}_{\bbQ_p}^{K^p}} = \lim_{K_p, n} \bbT'_{K_pK^p,n}.\] 
In particular, considering the maps from the limit when $n=1$ we deduce that the action on $\mc{A}^{K^p}_{\bbF_p}$ factors over $\bbT'^\wedge_{\mc{A}_{\bbQ_p}^{K^p}}$ as claimed. Moreover, each term in the limit has the discrete topology, so we also find that every \emph{open} maximal ideal is pulled back from $\bbT'_{K_pK^p,n}$ for some $n$, and we can take $n=1$ because any maximal ideal in $\bbT'_{K_pK^p,n}$ contains $p$. But the maximal ideals in $\bbT'_{K_pK^p,1}$ are precisely those such that $(\mc{A}_{\bbF_p}^{K_pK^p})_\frakm \neq 0$, so we conclude. 
\end{proof}

\begin{lemma}\label{lemma:hecke-action-open-p-adic-mf} The action of $\bbT'$ on $\mc{M}^{K^p}_{\bbF_p}$ factors through $\bbT'^\wedge_{\bbV^{K^p}_{\bbQ_p}}$ and the maximal ideals such that $\left(\mc{M}^{K^p}_{\bbF_p}\right)_\frakm \neq 0$ are precisely the open ideals of $\bbT'^\wedge_{\bbV^{K^p}_{\bbQ_p}}.$
\end{lemma}
\begin{proof}
Write $F_i \bbV_{\bbQ_p}^{K^p}$ for the image of $\bigoplus_{k=0}^{i} M_{k,\bbQ_p}^{K^p}$ in $\bbV^{K^p}_{\bbQ_{p}}$ under the evaluation map. The subspaces $F_i \bbV_{\bbQ_p}^{K^p}$ are preserved by $\bbT'$ and their union is dense in $\bbV_{\bbQ_p}^{K^p}$. If we write $F_i \bbV_{\bbZ_p}^{K^p} = F_i \bbV_{\bbQ_p}^{K^p} \cap \bbV_{\bbZ_p}^{K^p}$ and $\bbT'_{i,n}$ for the image of $\bbT'$ in $\End_{\bbZ/p^n} (F_i \bbV_{\bbZ_p}^{K^p} / p^n F_i \bbV_{\bbZ_p}^{K^p})$, then Lemma \ref{lemma:dense-prodiscrete} gives $\bbT'^\wedge_{\bbV_{\bbQ_p}^{K^p}} = \lim_{i,n} \bbT'_{i,n}.$

We write $F_i\bbV^{K^p}_{\bbF_p} := F_i \bbV_{\bbZ_p}^{K^p} / \left(F_i \bbV_{\bbZ_p}^{K^p}  \cap p\bbV_{\bbZ_p}^{K^p} \right) = F_i \bbV_{\bbZ_p}^{K^p} / p F_i \bbV_{\bbZ_p}^{K^p}$. Then, $F_i \bbV_{\bbF_p}^{K^p}$ is a finite dimensional $\bbF_p$-vector space and $\bbV_{\bbF_p}^{K^p} = \bigcup_i F_i \bbV_{\bbF_p}^{K^p}$. Taking the maps corresponding to the $n=1$ terms in the limit above, we deduce that the action of $\bbT'$ on $\bbV_{\bbF_p}^{K^p}$ factors through $\bbT'^\wedge_{\bbV_{\bbQ_p}^{K^p}}$. 

Now, because $\mc{M}^{K^p}_{\bbF_p}=(\bbV_{\bbF_p}^{K^p})^{1+p\bbZ_p}$, the action on $\mc{M}_{\bbF_p}^{K^p}$ also factors through $\bbT'^\wedge_{\bbV_{\bbQ_p}^{K^p}}.$  Moreover, arguing as in the proof of the previous lemma, any open maximal ideal is pulled back from $\bbT'_{i,1}$ for some $i$, i.e. from the image of $\bbT'$ in $\End(F_i \bbV_{\bbF_p}^{K^p})$, so we conclude that the open maximal ideals are exactly those for which $(F_i \bbV_{\bbF_p}^{K^p})_\frakm \neq 0$ for some $i$. On the other hand, $F_i \bbV_{\bbF_p}^{K^p}$ is preserved by the action of $\bbZ_p^\times$, which commutes with the Hecke action, thus $(F_i \bbV_{\bbF_p}^{K^p})_\frakm \neq 0$ is a finite dimensional $\bbF_p$-vector space with an action of $\bbZ_p^\times$. Because $1+p\bbZ_p$ is a pro-$p$ group, $(F_i \bbV_{\bbF_p}^{K^p})_\frakm \neq 0$ if and only if $(F_i \bbV_{\bbF_p}^{K^p})^{1+p\bbZ_p}_\frakm \neq 0$, as in the proof of Lemma \ref{lemma.d-to-np}. Because $\mc{M}_{\bbF_p}^{K^p} = \bigcup_i (F_i \bbV_{\bbF_p}^{K^p})^{1+p\bbZ_p}$, we find that as $i$ varies these are exactly the maximal ideals for which $(\mc{M}_{\bbF_p}^{K^p})_\frakm \neq 0.$ 
\end{proof}

Combining Lemmas \ref{lemma:hecke-action-open-dtimes} and \ref{lemma:hecke-action-open-p-adic-mf} (and invoking Lemma \ref{lemma.d-to-np} still to see that only finitely many maximal ideals appear), we obtain the $\bbF_p$-version of Theorem \ref{theorem:serre} as a consequence of Theorem \ref{mainthm:HeckeIso}.

\section{The spectral $p$-adic Jacquet-Langlands correspondence}\label{section:p-adic}

In this section we prove Theorem \ref{mainthm:HeckeIso}. After recalling some preliminaries on perfectoid modular curves  in  \S\ref{subsection:perfectoid-modular-curves}-\ref{subsection:hodge-tate-period-map}, in \S\ref{subsection:perfectoid-igusa-variety}-\ref{subsection:uniformization-of-fibers}  we recall the construction of the supersingular perfectoid Igusa variety and its relation to the fibers of the Hodge-Tate period map. In \S\ref{subsection:evaluation} we use this relation to construct an evaluation map from classical modular forms to $\calA_{\bbC_p}$. The key properties of this evaluation map are that it is injective and has dense image; these are established in Theorem \ref{theorem:eval-inj-dense}. We conclude in \S\ref{ss.completed-hecke-comparison} by invoking some of the tools developed in \S\ref{ss.completed-actions} -- in particular, Theorem \ref{theorem:list-of-isomorphisms} gives the isomorphism of completed Hecke algebras of Theorem \ref{mainthm:HeckeIso} along with some other isomorphisms indicated in the introduction.  

\subsection{Perfectoid modular curves}\label{subsection:perfectoid-modular-curves}
Let $K^p \leq \GL_2(\bbA_f^{(p)})$ be a sufficiently small compact open subgroup. Viewing $K^p$ as a subgroup of $\GL_2(\bbA_f)$, the (infinite level at $p$) schematic modular curves $Y_{K^p}$ and $X_{K^p}$ define functors on affinoid perfectoids over $(\mbb{C}_p, \mc{O}_{\bbC_p})$ by
\[ \mc{Y}_{K^p}: \Spa(R,R^+) \mapsto Y_{K^p}(R) \textrm{ and } \mc{X}_{K^p}: \Spa(R,R^+) \mapsto X_{K^p}(R). \]
It will be convenient to write out an equivalent definition for $\mc{Y}_{K^p}$, where we separate out the $p$ and prime-to-$p$ parts of the level:
\begin{definition}$\mc{Y}_{K^p}$ is the functor on affinoid perfectoid $\Spa(R,R^+)/\Spa(\bbC_p, \mc{O}_{\mbb{C}_p})$
\[ \mc{Y}_{K^p}: \Spa(R,R^+) \mapsto \{ (E, \varphi_p, \mc{K}^p) \} /\sim \]
sending $(R,R^+)$ to the set of equivalence classes of triples $(E, \varphi_p, \mc{K}^p)$ where
\begin{enumerate}
\item $E/R$ is an elliptic curve,
\item $\varphi_p: \ul{\bbQ_p}^2 \xrightarrow{\sim} V_{\bbQ_p}(E)$ is an isomorphism
\item $\mc{K}^p \subset \ul{\Isom}( (\ul{\bbA_f^{(p)}})^2, V_{\bbA_f^{(p)}}(E))$ is a $K^p$-torsor,
\item The relation $\sim$ is defined by $(E,\varphi_p, \mc{K}^p) \sim (E', \varphi_p', {\mc{K}^p}')$ if there is a quasi-isogeny $q: E \rightarrow E' $ such that $q \circ \varphi_p = \varphi_p'$ and $q(\mc{K}^p) = {\mc{K}^p}'.$
\end{enumerate}
\end{definition}

Both $\mc{Y}_{K^p}$ and $\mc{X}_{K^p}$ are represented by perfectoid spaces over $\mbb{C}_p$: in \cite[Theorem 3.3.18]{scholze:torsion} it is shown that there is a perfectoid space 
\[ \mc{X}_{K^p}' \sim \lim_{\substack{K_p \leq \GL_2(\mbb{Z}_p) \\ \textrm{compact open}}} X_{K_pK^p,\mbb{C}_p}^{\ad} \]
where here a $\sim$ limit, as defined in \cite[Definition 2.4.1]{scholze-weinstein:moduli-of-p-divisible-groups}, in particular implies (via \cite[Proposition 2.4.5]{scholze-weinstein:moduli-of-p-divisible-groups}) that for $(R,R^+)$ perfectoid over $(\mbb{C}_p, \mc{O}_{\mbb{C}_p})$, 
\[ \mc{X}_{K^p}'(R,R^+)= \lim_{\substack{K_p \leq \GL_2(\mbb{Z}_p) \\ \textrm{compact open}}} X_{K_pK^p,\mbb{C}_p}^{\ad}(R,R^+). \]
But we have the identity 
\[ X_{K_pK^p, \mbb{C}_p}^\ad(R, R^+)=X_{K_pK^p}(R) \]
because $X_{K_pK^p}$ is projective and the line bundles on $\Spa(R,R^+)$ and $\Spec R$ are identified via their global sections which are  rank one projectives over $R$ (for affine schemes this is standard, and for $\Spa(R,R^+)$ the equivalence between vector bundles and finite projective modules over $R$ is given by \cite[Theorem 8.2.2]{kedlaya-liu:foundations}). Thus we conclude $\mc{X}_{K^p}=\mc{X}_{K^p}',$ so that $\mc{X}_{K^p}$ is a perfectoid space.  

Passing to the open modular curves we find
\[ \lim_{\substack{K_p \leq \GL_2(\mbb{Z}_p) \\ \textrm{compact open}}} Y_{K_pK^p,\mbb{C}_p}^{\ad}. \]
is also represented by an open in $\mc{X}_{K^p}$, which is perfectoid by the above result. It is immediate by construction of the analytification of an affine finite type scheme over $\mbb{C}_p$ that 
\[ Y_{K_pK^p,\mbb{C}_p}^{\ad}(R,R^+)=Y_{K_pK^p}(R), \]
so that we conclude $\mc{Y}_{K^p}$ is a perfectoid space.

\begin{remark} For open modular curves we can then obtain the same result for $K^p$ closed but not necessarily compact open by using Proposition \ref{prop:modular-curves} and the fact that the limit of a tower of finite \'{e}tale covers of perfectoid spaces is perfectoid, and similarly for $X$ if this is combined with explicit computations at the cusps.  Since our main concern is the Hecke action, for which it suffices to have a construction for $K^p$ compact open, we do not go further into these details here. 
\end{remark} 

From the actions on $Y_{K^p}$ and $X_{K^p}$, we obtain an action of $\GL_2(\mbb{Q}_p)$ on $\mc{X}_{K^p}$ preserving $\mc{Y}_{K^p}$ and a commuting action of $\GL_2(\bbA_f^{(p)})$ on the towers as $K^p$ varies over sufficiently small compact opens  $K^p \leq \GL_2(\bbA_f^{(p)})$.

\subsection{Modular forms}\label{subsection.modular-forms-last-section}

For $k \geq 0$ recall from \S\ref{ss.modular-forms} that we have the admissible $\GL_2(\bbA_f)$-representation of weight $k$ modular forms over $\bbC_p$,
\[ M_{k,\bbC_p} = H^0(X_{\mbb{C}_p}, \omega^k)\cong\colim_{\substack{K \leq \GL_2(\bbA_f) \\ K \textrm{ sufficiently small} \\\textrm{compact open}}} H^0(X_{K,\bbC_p}, \omega^k). \]
We now relate these to $\mc{X}_{K^p}$. 

The $\GL_2(\mbb{Q}_p)$-action on $\mc{X}_{K^p}$ can be viewed as induced from the identification 
\[ \mc{X}_{K^p} \sim \lim_{\substack{K_p \leq \GL_2(\mbb{Q}_p) \\ \textrm{compact open}}} X^\ad_{K_pK^p, \mbb{C}_p}  \] 
and the action of $\GL_2(\mbb{Q}_p)$ on the tower. In particular, the analytification of the $\GL_2(\mbb{Q}_p)$-equivariant modular bundle $\omega$ 
on the tower pulls back to a $\GL_2(\mbb{Q}_p)$-equivariant vector bundle on $\mc{X}_{K^p}$ which we also denote by $\omega$. Moreover, these compile to a $\GL_2(\bbA_f^{(p)})$-equivariant bundle on the tower $(\mc{X}_{K^p})_{K^p}$ as $K^p$ varies. 

For $K_p \leq \GL_2(\mbb{Q}_p)$ compact open we have 
\[ H^0(X_{K_pK^p, \mbb{C}_p}^\ad, \omega^k)=H^0(X_{K_pK^p, \mbb{C}_p}, \omega^k) = M_{k, \bbC_p}^{K_pK^p}. \]
In particular, by pullback to $\mc{X}_{K^p}$ we obtain a $\GL_2(\mbb{Q}_p)$-equivariant injection 
\[ M_{k, \bbC_p}^{K^p} = H^0(X_{K^p}, \omega^k) \hookrightarrow H^0(\mc{X}_{K^p}, \omega^k ). \]
Though we will not need any more than this, we note that this is an isomorphism onto the $\GL_2(\mbb{Q}_p)$-smooth vectors, as can be checked after restriction to $\mc{Y}_{K^p}$ where it follows from the sheaf property for the completed structure sheaf on the pro-\'{e}tale site of $\mc{Y}_{K_pK^p, \bbC_p}^\ad$ for each compact open subgroup $K_p \leq \GL_2(\mbb{Q}_p)$. 

This isomorphism is compatible with the $\GL_2(\bbA_f^{(p)})$-action on the tower of $\mc{X}_{K^p}$ so extends to a $\GL_2(\bbA_f)$-equivariant injection 
\[ M_{k, \mbb{C}_p} \hookrightarrow \colim_{K^p} H^0(\mc{X}_{K^p}, \omega^k). \]

\begin{remark} As in Remark \ref{remark:mod-p-twist}, to obtain a $\GL_2(\bbA_f)$-action inducing the standard Hecke action, one should introduce a twist by the unramified determinant character that appears in the Kodaira-Spencer isomorphism. This does not change the resulting completed Hecke algebra, so it is reasonable to include the twist only when convenient, that is, when comparing with singular/\'{e}tale cohomology.  
\end{remark}

\subsection{The Hodge-Tate period map}\label{subsection:hodge-tate-period-map}

In \cite{scholze:torsion}, Scholze constructs the Hodge-Tate period map, 
\[ \pi_{\HT, K^p}: \mc{X}_{K^p} \rightarrow \PP^1 (= \bbP^{1,\ad}_{\bbC_p}). \]
We normalize some choices related to the group actions and equivariant structures by requiring that, over $\mc{Y}_{K^p}$, $\pi_{\HT,K^p}$ is the classifying map for the line 
\[ \Lie E(1) \subset V_p E \otimes \Oc_{\mc{Y}_{K^p}} \isoeq \Oc_{\mc{Y}_{K^p}}^2 \]
given by the Hodge-Tate filtration for the universal elliptic curve and the trivialization of its relative Tate module. In particular, $\pi_{\HT,K^p}$ is $\GL_2(\mbb{Q}_p)$-equivariant for the action on $\bbP^1$ in which $\GL_2(\bbQ_p)$ acts by the \emph{dual} of the standard representation on $H^0(\bbP^1, \mc{O}(1))$. The image of the boundary/cusps $\mc{X}_{K^p}\bs \mc{Y}_{K^p}$ is $\mbb{P}^1(\mbb{Q}_p)$.

The maps $\pi_{\HT,K^p}$ compile to a $\GL_2(\bbA_f)$-equivariant map 
\[ \pi_\HT: (\mc{X}_{K^p})_{K^p} \rightarrow \PP^1 \]
from the tower $\mc{X}_{K^p}$, where $\PP^1$ is given the trivial action of $\GL_2(\bbA_f^{(p)})$. By construction, after fixing a compatible system of $p$-power roots of unity in $\bbC_p$, there is a natural isomorphism $\pi_{\HT}^{*} \calO(1) = \omega$ of $\GL_2(\AA_f)$-equivariant line bundles on the tower $(\mc{X}_{K^p})_{K^p}$.

\begin{remark} Comparing with \cite{caraiani-scholze:generic, scholze-weinstein:moduli-of-p-divisible-groups}, over the good reduction locus one can also think of the map $\pi_\HT$  as the classifying map for the $p$-divisible group $E[p^\infty]$ equipped with a basis for $T_p E[p^\infty].$ A key property of $\pi_\HT$ is that we can construct $\GL_2(\bbA_f^{(p)})$-equivariant fake Hasse invariants via pullback; this matches well with the perspective on the classical Hasse invariant described in \ref{sss.hasse} which can be read  as saying that it is pulled back from the moduli stack of $p$-divisible groups. 
\end{remark}

\subsection{The perfectoid supersingular Igusa variety}\label{subsection:perfectoid-igusa-variety}
\newcommand{\Igfoid}{\mc{I}g}
For a closed subgroup $K^p \leq \GL_2(\mbb{A}_f^{(p)})$, the supersingular Caraiani-Scholze Igusa moduli problem of Definition \ref{def.igusa} defines a functor on affinoid perfectoids over $(\mbb{C}_p, \mc{O}_{\mbb{C}_p})$,
\[ \Igfoid^\ssing_{K^p}: \Spa(R,R^+) \rightarrow \Ig_{K^p}^\ssing(R^+/p). \]
When $K^p$ is sufficiently small, Corollary \ref{corollary.igusa-level-coset} identifies $\Ig_{K^p}^\ssing$ with $\Spec A$ for 
\[ A=\mc{A}_{\overline{\bbF}_p}^{K^p} = \Cont(D^\times(\mbb{Q})\bs D^\times(\bbA_f)/K^p, \overline{\bbF}_p). \]
Because $A$ is  a perfect ring and $R^+$ is $p$-adically complete,
\begin{align*} \Ig_{K^p}^\ssing(R^+/p) &= \Hom_{\overline{\bbF}_p}(A, R^+/p) \\ &= \Hom_{\breve{\bbZ}_p}(W(A), R^+) \\ &= \Hom_{(\breve{\bbQ}_p, \breve{\bbZ}_p)}\left((W(A)[1/p], W(A)), (R,R^+) \right) \\
&= \Hom_{(\breve{\bbQ}_p, \breve{\bbZ}_p)}\left((\mc{A}_{\breve{\bbQ}_p}^{K^p}, \mc{A}_{\breve{\bbZ}_p}^{K^p}), (R,R^+) \right) \\
&= \Hom_{(\bbC_p, \mc{O}_{\bbC_p})}\left((\mc{A}_{\bbC_p}^{K^p}, \mc{A}_{\mc{O}_{\mbb{C}_p}}^{K^p}), (R,R^+) \right). 
\end{align*}
But $(\mc{A}_{\bbC_p}^{K^p}, \mc{A}_{\mc{O}_{\mbb{C}_p}}^{K^p})$ is a perfectoid Huber pair, so we conclude that $\Igfoid^\ssing_{K^p}$ is represented by the affinoid perfectoid space $\Spa(\mc{A}_{\bbC_p}^{K^p}, \mc{A}_{\mc{O}_{\mbb{C}_p}}^{K^p})$, which is just the profinite set $D^\times(\mbb{Q})\bs D^\times(\bbA_f)/K^p$ viewed as a perfectoid space over $\Spa(\mbb{C}_p, \mc{O}_{\bbC_p})$.

\begin{remark} In \cite[Section 4]{caraiani-scholze:generic} this procedure is carried out in much greater generality to construct perfectoid Igusa varieties (without the explicit identification with a double coset, which in general will only occur over the basic locus). 
\end{remark}

\subsection{Uniformization of fibers of $\pi_{\HT}$} \label{subsection:uniformization-of-fibers}
\newcommand{\loc}{\mathrm{loc}}
We fix now a one-dimensional connected height two  $p$-divisible group $G / \calO_{\bbC_p}$ equipped with a quasi-isogeny
\[ \rho_{G}: \bbX_{\ssing, \calO_{\bbC_p}/p} \otimes \bbQ_p \xrightarrow{\sim} G_{\calO_{\bbC_p}/p} \otimes \bbQ_p \] 
and a trivialization $\psi_{G}:  \bbZ_p^2 \xrightarrow{\sim} T_p G (\calO_{\bbC_p})$. By the Scholze-Weinstein classification \cite{scholze-weinstein:moduli-of-p-divisible-groups}, the pair $(G,\psi_G)$ is equivalent to the point $x \in \bbP^1(\bbC_p) \bs \bbP^1(\bbQ_p)$ determined by the position of the Hodge-Tate filtration. As it will be convenient later, we choose our data so that $x$ is in the affinoid ball $B_1: |U/V| \leq 1$, where here $U$ and $V$ denote the standard basis for $O(1)$ (i.e. projective coordinates are $[U:V]$). 

As in \cite{scholze-weinstein:moduli-of-p-divisible-groups}, there is a perfectoid Lubin-Tate space $\widehat{\mc{M}}_{\LT,\infty}$ over $\Spa(\mbb{C}_p, \mc{O}_{\mbb{C}_p})$ parameterizing quasi-isogeny lifts $G$ of $\bbX_\ssing$ equipped with trivializations of $T_pG$. We do not define this space more carefully as we use it only through a citation to \cite{caraiani-scholze:generic} in a proof below, but let us recall that there is a local Hodge-Tate period map classifying the Hodge-Tate filtration,
\[ \pi_{\HTloc}: \widehat{\mc{M}}_{\LT,\infty} \rightarrow \mbb{P}^1, \]
and $\rho_G$ determines a point $x_\infty \in \widehat{\mc{M}}_{\LT,\infty}(\bbC_p)$ lying in $\pi_{\HTloc}^{-1}(x)$. 

From this data we now construct a map
\[ \unif_{x_\infty, K^p}: \Igfoid_{K^p} \rightarrow \mc{Y}_{K^p} \subset \mc{X}_{K^p}.\]
On points in affinoid perfectoid $(R,R^+)$, it sends $(E, \varphi_p, \mc{K}^p)$ to $(E', \varphi_p', {\mc{K}^p}')$, where
\begin{enumerate}
\item $E'=\mc{E}_R$ where $\mc{E}$ is the elliptic curve over $R^+$ given by the Serre-Tate lifting (cf. \cite[Section 1]{katz:serre-tate-local-moduli}) of $E$ to $R^+$ determined by the quasi-isogeny
\[ \rho_G \circ \varphi_{p}^{-1}: E[p^\infty] \otimes \bbQ_p \rightarrow G_{R^+/p} \otimes \bbQ_p. \] 
\item $\varphi_p'$ is the composition of $\psi_G$ with the canonical identification of Tate modules induced by $\mc{E}[p^\infty]=G[p^\infty]_{R^+}$. 
\item ${\mc{K}^p}'$ is the unique lift of $\mc{K}^p$. 
\end{enumerate}
By construction, $\unif_{x_\infty}$ factors through the closed subset $\pi_{\HT, K^p}^{-1}(x) \subset \mc{X}_{K^p, \bbC_p}$. The latter is a Zariski closed subset of a perfectoid space, and thus admits a canonical structure of a perfectoid space through which $\unif_{x_\infty}$ factors (because the domain is also a perfectoid space). 

In addition to the obvious $\GL_2(\bbA_f^{(p)})$-equivariance in the tower, the map $\unif_{x_\infty,K^p}$ also satisfies an equivariance at $p$: let $T_x:=\Aut(G \otimes \bbQ_p)$ be the group of self quasi-isogenies of $G$. Then, by the Scholze-Weinstein classification, $\psi_G$ identifies $T_x$ with the stabilizer in $\GL_2(\mbb{Q}_p)$ of $x\in \bbP^1(\mbb{C}_p) \backslash \bbP^1(\mbb{Q}_p)$. On the other hand, $\rho_G$ identifies $T_x$ with a subgroup of the self quasi-isogenies
\[ \Aut(\mbb{X}_{\ssing, \mc{O}_{\mbb{C}_p}/p} \otimes \bbQ_p) = \Aut(\mbb{X}_{\ssing, \overline{\bbF}_p} \otimes \mbb{Q}_p) = D_p^\times. \] 
By the interpretation as a stabilizer, $T_x$ is equal to either $\QQ_p^\times$ or $F^\times$ for a quadratic extension $F/\bbQ_p$, and the latter occurs exactly when $x \in \PP^1(F) - \PP^1(\QQ_p)$ -- this follows from an explicit computation after observing that any line preserved by a non-scalar matrix in $\GL_2(\mbb{Q}_p)$ must be defined over a quadratic extension of $\bbQ_p$. 

\begin{theorem}\label{theorem.perfectoid-isomorphism-fiber}
The maps $\unif_{x_\infty, K^p}$ induce a $T_x\times\GL_2(\bbA_f^{(p)})$-equivariant isomorphism of towers of perfectoid spaces
\[ \unif_{x_\infty}: \left( D^\times(\bbQ) \bs D^\times(\bbA_f)/K^p \right)_{K^p} = \left(\Igfoid^\ssing_{K^p}\right)_{K^p} \xrightarrow{\sim} \left(\pi_{\HT,K^p}^{-1}(x)\right)_{K^p} \]
as $K^p$ varies over sufficiently small compact open subgroups of $\GL_2(\bbA_f^{(p)})$. 
\end{theorem}
\begin{proof}
The equivariance and compatibility in the tower is immediate from the definition, so it suffices to show that for each $K^p$ the map $\unif_{x_\infty, K^p}$ is an isomorphism of perfectoid spaces. 

We first note that $\unif_{x_\infty,K^p}$ is the restriction to $\Igfoid^\ssing_{K^p} \times x_\infty$ of a map
\[ \unif_{K^p}:  \Igfoid^\ssing_{K^p} \times_{\Spa(\CC_p,\Oc_{\CC_p})} \widehat{\mc{M}}_{\LT,\infty} \rightarrow \mc{Y}_{K^p}^\ssing \]
defined similarly in \cite[Section 4]{caraiani-scholze:generic} --  here $\mc{Y}_{K^p}^\ssing$ denotes the closed subspace of $\mc{Y}_{K^p}$ whose $\Spa(C,C^+)$ points for $C$ complete algebraically closed correspond to triples $(E, \rho_p, \mc{K}^p)$ such that $E$ extends to an elliptic curve over $C^+$ with supersingular reduction at the maximal ideal ${\frakm}_{C^+}$ of $C^+$. 

 It follows from \cite[Lemma 4.3.20]{caraiani-scholze:generic} (cf. also \cite[Definition 4.3.17]{caraiani-scholze:generic}), that the diagram
\[ \xymatrix{ \Igfoid^\ssing_{K^p} \times_{\Spa(\CC_p,\Oc_{\CC_p})} \widehat{\mc{M}}_{\LT,\infty} \ar[r] \ar[d]^{\unif_{K^p}} & \widehat{\mc{M}}_{\LT,\infty} \ar[d]^{\pi_\HTloc} \\
\mc{Y}_{K^p}^\ssing \ar[r]^{\pi_{\HT,K^p}} & \PP^1 
} \]
is Cartesian on perfectoid spaces,  Thus, because $\pi_\HTloc(x_\infty)=x$, we find that $\unif_{x_\infty, K^p}$ induces an isomorphism between the functors of points of $\Igfoid^\ssing_{K^p}$ and $(\pi_{\HT,K^p}|_{\mc{Y}_{K^p}^\ssing})^{-1}(x)$ on perfectoid spaces. Because both are perfectoid spaces, and a perfectoid space is determined by its functor of points on perfectoid spaces, we conclude that $\unif_{x_\infty}$ is an isomorphism onto this space.

It remains to show that $\pi_{\HT,K^p}^{-1}(x) \subset \mc{Y}_{K^p}^\ssing$. Here the point is just that the loci of ordinary and multiplicative reduction map to $\bbP^1(\bbQ_p)$, but let us make this completely precise: suppose given a geometric point 
\[ (E, \varphi_p, \mc{K}^p) \in \pi_{\HT,K^p}^{-1}(x)(C,C^+), \]
and write $t$ for the corresponding point 
\[ t=(E,\mc{K}^p) \in X_{\GL_2(\bbZ_p)K^p}(C)=\Xint_{K^p}(C^+). \]
By \cite[Lemma 3.3.19]{scholze:torsion}, because $x \not\in \bbP^1(\bbQ_p)$, there is a rational number $1\geq \epsilon>0$ such that the Hasse invariant vanishes on $t_{C^+/p^\epsilon}$, which thus factors through the supersingular locus. But this implies in particular that $E$ has  supersingular reduction over $C^+/\frakm_{C^+}$, so we conclude. 

\end{proof}

\subsection{Evaluation of modular forms}\label{subsection:evaluation}
Any choice of basis for the fiber of $\calO(1)$ at $x \in \bbP^1(\bbC_p)$ induces a $\GL_2(\bbA_f^{(p)})$-equivariant trivialization of $\omega|_{\pi_{\HT}^{-1}(x)}$. In particular, if we take as basis the section $V|_x$ (recall $U$ and $V$ are the standard basis for $O(1)$), we obtain a map from classical modular forms
\[ \eval_{x_\infty}^{K^p}:  \bigoplus_{k \geq 0} M_{k,\bbC_p}^{K^p} \rightarrow \calA_{\bbC_p}^{K^p} \l( = H^0(\Igfoid_{K^p}^\ssing, \calO) \r) \]
by first pulling back to $\mc{X}_{K^p}$ then restricting to $\pi_{\HT}^{-1}(x)$ and dividing by $V^k$. 

\begin{remark} Choosing $G$ to have endomorphisms by $\breve{\bbZ}_{p^2}$ and scaling the trivialization $V|_x$ by a $p$-adic period of $G$, one can obtain an evaluation map that, on $H^0(\Xint_{\breve{\bbZ}_p, \omega^k}) \subset M_{k,\bbC_p}^{K^p}$, reduces modulo $p$ to Serre's evaluation map described in \S\ref{section:mod-p}. \end{remark}

\begin{theorem}\label{theorem:eval-inj-dense}
For $K^p \leq \GL_2(\bbA_f^{(p)})$ a sufficiently small compact open, the map $\eval_{x_\infty}^{K^p}$ is injective and has dense image. Moreover, as $K^p$ varies they are compatible and induce an injective
$\GL_2(\bbA_f^{(p)})$-equivariant map 
\[ \eval_{x_\infty}: \bigoplus_{k\geq0}M_{k,\bbC_p} \rightarrow \mc{A}_{\bbC_p} \]
\end{theorem}
\begin{proof}
The compatibility and equivariance as $K^p$ varies is clear from construction.

We now show injectivity: because $z \in \bbZ_p^\times \subset T_x$ acts as multiplication by $z^{-k}$ on our fixed trivialization $V|_x$, we find that for any compact open $K_p \subset \GL_2(\bbQ_p)$, the image of $M^{K_pK^p}_{k, \bbC_p}$ in $H^0(\pi_\HT^{-1}(x), \mc{O})$ transforms under $z^{k}$ for the action of $z \in K_p \cap \bbZ_p^\times$. Thus, by the equivariance of the $T_x$ action in Theorem \ref{theorem.perfectoid-isomorphism-fiber}, the image of $M^{K^p}_{k, \bbC_p}$ lands in the subspace $\calA^{K^p}_{\bbC_p}[k]$ of vectors on which the action of the central $\bbZ_p^\times \leq D_p^\times$ is differentiable with derivative $k$. In particular, there can be no cancellation between the different degrees $k$.

To show the map is injective on each $M_{k,\bbC_p}$, we first observe that for $K=K_pK^p$ a compact open subgroup of $\GL_2(\bbA_f)$, the image of $\pi_{\HT}^{-1}(x)$ in $X_{K, \bbC_p}^\ad(\bbC_p)=X_{K, \bbC_p}(\bbC_p)$ intersects every connected component of $X_{K,\bbC_p}$ in an infinite set -- indeed, the map factors as an injection from $D^\times(\bbQ)\bs D^\times(\bbA_f) / (K_p \cap T_x) \cdot K^p$, and the connected component of the image in $X_{K,\bbC_p}$ is recorded by the map
\[ g = g_p \times (g_\ell)_{\ell \neq p} \mapsto \Nrd g_p \cdot |\Nrd g_p|_p \cdot \prod_{\ell \neq p} |\det g_\ell|_\ell \]
with values in $\bbZ_p^\times$ modulo the image of $K$. 

If $s$ is a non-zero section of a line bundle on $X_{K,\bbC_p}$ then there is at least one connected component where it has only finitely many zeroes. Thus, any section of $\omega^k$ over $X_{K}$ which vanishes upon restriction to $\pi_{\HT,K^p}^{-1}(x)$ is identically zero, and we conclude the map is injective. 

We now show the map has dense image. By assumption, $x$ is contained in
\[ B_1 :|U|/|V| \leq 1 \subset \bbP^1, \]
and \cite[Theorem 3.3.18(i)]{scholze:torsion} gives that $\pi_{\HT}^{-1}(B_1)$ is affinoid perfectoid. By \cite[Remark 7.5]{bhatt-scholze:prisms} (Zariski closed implies strongly Zariski closed; cf. also \cite[Definition II.2.6]{scholze:torsion}), we find that the map 
\[ H^0(\pi_{\HT,K^p}^{-1}(B_1), \calO^+) \rightarrow H^0(\pi_{\HT,K^p}^{-1}(x), \calO^+) = \calA_{\calO_{\bbC_p}}^{K^p} \]
is almost surjective; in particular, the image contains $p \calA_{\calO_{\bbC_p}}^{K^p}$.

Thus, given ${f \in p \cdot \calA_{\calO_{\bbC_p}}^{K^p}}$, we can lift it to $\tilde{f} \in H^0(\pi_{\HT,K^p}^{-1}(B_1), \calO^+) $. By Lemma \ref{lem.fake-Hasse-lifting} below we find that for any $n > 0$ that there is a $k$ large enough, $K_p$ small enough, and an element of $\alpha \in M_{k,\bbC_p}^{K_pK^p}$ such that 
\[ \alpha/V^k \in H^0(\pi_{\HT,K^p}^{-1}(B_1), \calO^+),\; \alpha/V^k \equiv \tilde{f} \textrm{ mod } p^n H^0(\pi_{\HT,K^p}^{-1}(B_1),\calO^+). \]
In particular, $\eval_{x_\infty}^{K^p}(\alpha) \equiv f \mod p^n$. Thus the image of $\eval_{x_\infty}^{K^p}$ contains a dense subset of the open ball $p\mc{A}^{K^p}_{\mc{O}_{\bbC_p}}$, so is dense in $\mc{A}^{K^p}_{\mbb{C}_p}$. 

\end{proof}

\begin{corollary}\label{corollary.evaluatoin-any-level} If $K^p \leq \GL_2(\bbA_f^{(p)})$	is \emph{any} compact open subgroup, then $\eval_{x_\infty}$ restricts to a Hecke-equivariant injection with dense image 
\[ \bigoplus_{k \geq 0} M_{k, \bbC_p}^{K^p} \rightarrow \mc{A}_{\bbC_p}^{K^p}. \] 
\end{corollary}
\begin{proof}
This statement is immediate from Theorem \ref{theorem:eval-inj-dense} except for the density when $K^p$ is not sufficiently small. To see this density, we pass to a sufficiently small compact open  normal subgroup ${K^p}' \trianglelefteq K^p$, then average approximations by ${K^p}'$-invariants over $K^p/{K^p}'$ to obtain approximations by $K^p$ invariants. 
\end{proof}

The following lemma extracts part of the fake Hasse invariant argument given in \cite[Proof of Theorem 4.3.1, pp. 1028-1031]{scholze:torsion} in the simplest possible case. We make no new contribution, however, as the specific statement we need is significantly simpler than the setup in loc. cit., we reproduce the proof here in the hopes that it will be helpful for the reader.   
\begin{lemma}\label{lem.fake-Hasse-lifting}	
If $g \in H^0(\pi_{\HT, K^p}^{-1}(B_1), \mc{O}^+)$ and $n>0$, then there is a compact open subgroup $K_p \subset \GL_2(\bbZ_p)$, a $k>0$, and an $\alpha \in H^0(X_{K_pK^p}, \omega^k)$ such that $\alpha/V^k \in H^0(\pi_{\HT, K^p}^{-1}(B_1), \mc{O}^+)$ and $\alpha/V^k \equiv g \mod p^n.$ 
\end{lemma}
\begin{proof}
We will use density of sections to approximate $g$, $U$ and $V$ at finite level, then make an argument using formal models and ampleness of $\omega$ to extend the approximation of $g$ to a global section of a power of $\omega$ after multiplying by a large enough power of the approximation of $V$. 

To better match the notation of loc. cit., we write
\[ s_1=\pi_{HT,K_p}^*U,\, s_2=\pi_{\HT,K^p}^*V,\, \mc{U}_1 = |s_2/s_1| \leq 1 \subset \mc{X}_{K^p},\textrm{ and } \mc{U}_2=|s_1/s_2|\leq 1 \subset \mc{X}_{K^p}. \]
In particular, $\mc{U}_2=\pi_{\HT,K^p}^{-1}(B_1)$ in the previous notation. 

By \cite[Theorem 3.3.18(i)]{scholze:torsion}, both $\mc{U}_1$ and $\mc{U}_2$ are affinoid perfectoid and come from finite level affinoids $\mc{U}_{1,K_p}$ and $\mc{U}_{2,K_p}$ in $X_{K_pK^p,\bbC_p}^\ad$ for small enough compact open $K_p \leq \GL_2(\bbZ_p)$. Moreover they satisfy
\[ \colim_{K_p} \mc{O}(\mc{U}_{\bullet,K_p})  \subset \mc{O}(\mc{U}_{\bullet}) \textrm{ is a dense subset for }\bullet=1,2.\]
The same density statement holds for sections of $\omega$, thus we can choose, for $K_p$ small enough, sections 
\[ s_{j}^{(i)} \in H^0(\mc{U}_{i,K_p}, \omega),\; i,j=1,2 \]
such that $s_{i}^{(i)}$ is non-vanishing and the sections $s_{j}^{(i)}$ satisfy
\[ \left|\frac{s_j - s_{j}^{(i)}}{s_{i}^{(i)}}\right| \leq |p^n| \textrm{ on } \mc{U}_i. \]
By possibly taking $K_p$ smaller still, we may also assume that there is a function 
\[ h \in H^0(\mc{U}_{2,K_p}, \mc{O}^+) \textrm{ such that } |h - g| \leq |p^n| \textrm{ on } \mc{U}_2, \]
where we recall that $g$ is the function we are ultimately trying to approximate. 

\newcommand{\strange}{\mathrm{strange}}
Now, \cite[Lemma 2.1.1]{scholze:torsion} gives:
\begin{enumerate}
	\item a formal model $\frak{X}^\strange$ over $\Spf \mc{O}_{\bbC_p}$ for $X_{K_pK^p, \bbC_p}^\ad,$ 
	\item a cover of $\frak{X}^\strange$ by affines $\frak{U}_\bullet=\Spf \mc{O}^+(\mc{U}_{\bullet, K_p})$, $\bullet=1,2$, and
	\item an ample line bundle $\frak{w}/\frak{X}^\strange$ modeling $\omega/X_{K_pK^p, \bbC_p}^\ad$ and such that 
	 \[ \frak{w}(\frak{U}_i) = \mc{O}^+(\mc{U}_{\bullet, K_p}) \cdot s_{i}^{(i)}. \]
\end{enumerate}
In particular, the sections $s_{i}^{(j)}$ glue mod $p^n$ to a section $\overline{s_i}$ of $\frak{w}/p^n$, and $\frak{U}_i$ is exactly the locus where $\overline{s_i}$ is invertible.  We deduce by ampleness that for large enough $k$:
\begin{enumerate}
\item $H^1(\frak{X}^\strange, \frak{w}^k/p^n)=0$ and
\item $\overline{h} \cdot \overline{s_{2}}^k$ extends to a global section of $H^0(\frak{X}^\strange, \frak{w}^k/p^n).$ 
\end{enumerate}
The first item implies that $H^0(\frak{X}^\strange, \frak{w}^k) \rightarrow H^0(\frak{X}^\strange, \frak{w}^k/p^n)$ is surjective: consider the exact sequence of cohomology coming from 
\[ 0 \rightarrow \frak{w}^k/p^n \rightarrow \frak{w}^k/p^{2n} \rightarrow \frak{w}^k/p^n \rightarrow 0\]
to lift to $H^0(\frak{X}^\strange, \frak{w}^k/p^{2n})$ and get $H^1(\frak{X}, \frak{w}^{k}/p^{2n})=0$, then repeat and pass to the limit. Thus we can lift $\overline{h}\cdot \overline{s_2}^k$ to 
\[ \alpha \in H^0(\frak{X}^\strange, \frak{w}^k) \subset H^0(X_{K_pK^p,\bbC_p}^\ad, \omega^k),\]
 and this is the desired section. 
\end{proof}

\subsection{Comparison of completed Hecke algebras}\label{ss.completed-hecke-comparison}

We now prove Theorem \ref{mainthm:HeckeIso}: Let $K^p \leq \GL_2(\bbA_f^{(p)})= D^\times(\bbA_f^{(p)})$ be compact open, let 
\[ \bbT^\abs_{K^p} = \bbZ [K^p \bs \GL_2(\bbA_f^{(p)}) / K^p ], \]
and let $\bbT' \subset \bbT^\abs_{K^p}$ be a subring. Corollary \ref{corollary.evaluatoin-any-level} and Lemma \ref{lem:CompDense} give 
\[ \bbT'^\wedge_{\l(M_{k,\bbC_p}^{K_p K^p}\r)_{K_p, k}} = \bbT'^\wedge_{\calA_{\bbC_p}^{K^p}}.\]
 Invoking Lemma \ref{lem:CompBaseChange}, we find that we can replace $\bbC_p$ with $\bbQ_p$ on both sides.

On the other hand, combining  Lemmas \ref{lemma:density-p-adic-mf} and \ref{lem:CompDense}, we find that
\[ \bbT'^\wedge_{\l(M_{k,\bbQ_p}^{\GL_2(\bbZ_p)K^p}\r)_{k}} = \bbT'^\wedge_{\bbV^{K^p}_{\bbQ_p}}. \]
Thus, to deduce Theorem \ref{mainthm:HeckeIso}, it remains only to show 
\[ \bbT'^\wedge_{\l(M_{k,\bbQ_p}^{\GL_2(\bbZ_p)K^p}\r)_{k}} = \bbT'^\wedge_{\l(M_{k,\bbQ_p}^{K_p K^p}\r)_{K_p, k}}. \]
This is basically a well-known result of Hida  \cite[Equation (1.7)]{hida:galois}, however, we are not aware of a full proof in the literature. Thus we include a proof here, arguing with completed cohomology as explained by Emerton \cite[Remarks 5.4.2 and 5.4.3]{emerton:local-global-proof}. This argument forms part of the proof of the following result, which encodes all of the isomorphisms indicated above and in the introduction:

\begin{theorem}\label{theorem:list-of-isomorphisms}
The identity map $\bbT'\rightarrow \bbT'$ extends to a topological isomorphism of the completed Hecke algebras $\bbT'$ acting on the following:
\begin{enumerate}
\item $\l( M_{k,\bbQ_p}^{K_p K^p} \r)_{k, K_p}$ for $k$ varying over all non-negative integers and $K_p$ varying over all compact open subgroups of $\GL_2(\bbQ_p)$,
\item $\l( M_{k,\bbQ_p}^{\GL_2(\bbZ_p)K^p} \r)_k $ for $k$ varying over all non-negative integers, 
\item $\l( M_{2,\bbQ_p}^{K_pK^p} \r)$, for $K_p$ varying over all compact open subgroups of $\GL_2(\bbQ_p)$,
\item the completed cohomology of the modular curve at level $K^p$ (cf. \cite{emerton:local-global-conjecture}), 
\[ \widehat{H}^1_{K^p} = \tilde{H}^1_{K^p}:= \left(\varprojlim_m \varinjlim_{K_p} H^i(Y_{K_p K^p}(\CC), \ZZ/p^m)\right)[1/p], \]  
\item the space of quaternionic automorphic forms $\calA_{\bbQ_p}^{K^p}$, and
\item  the space of $p$-adic modular forms $\bbV^{K^p}_{\bbQ_p}$.
\end{enumerate}
\end{theorem}
\begin{proof}
We have already established the identities of completed Hecke algebras (1)=(5) and (2)=(6). We will conclude by identifying (4)=(1), (4)=(2), and (4)=(3). 

\newcommand{\aux}{\mathrm{aux}} 
For $K=K_pK^p$ sufficiently small, let $\mathrm{pr}:E_K \rightarrow Y_K(\CC)$ be the universal elliptic curve over $Y_K(\CC)$ and let $\underline{\Sym^k}$ be the $k$th symmetric power of $R^1 \mathrm{pr}_* \bbQ$, a $\GL_2(\AA_f)$-equivariant local system on the tower $(Y_K(\CC))_K$. Fixing an isomorphism $\overline{\QQ}_p \isoeq \CC$,
we obtain for each $k\geq 2$ and $K_p$ such that $K=K_pK^p$ is sufficiently small a $\TT^\abs$-equivariant injection
\begin{equation}\label{equation:funny-injection-completed-hecke} M^{K_pK^p}_{k, \QQ_p} \hookrightarrow H^1(Y_{K_pK^p}(\CC), \underline{\Sym^{k-2}})\otimes \CC_p \end{equation}
given by composing the maps
\begin{multline*}
M^{K_pK^p}_{k, \QQ_p} \hookrightarrow M^{K_p K^p}_{k, \CC} \hookrightarrow H^1(Y_{K_pK^p}(\CC), \underline{\Sym^{k-2}})\otimes \CC  \\ \hookrightarrow H^1(Y_{K_pK^p}(\CC), \underline{\Sym^{k-2}})\otimes \CC_p,
\end{multline*}
where the second arrow comes from the classical Eichler-Shimura isomorphism and the last arrow comes from composition of the isomorphism $\overline{\QQ_p}\isoeq \CC$ with $\overline{\QQ_p} \hookrightarrow \CC_p$. It follows from the Eichler-Shimura isomorphism that (\ref{equation:funny-injection-completed-hecke}) induces an isomorphism on the image of $\TT'$ in the respective endomorphism rings. If we denote $\TT'^\wedge_\aux$ the completed Hecke algebra for $\TT'$ acting on the family
\[ \left( H^1(Y_{K_p K^p}(\CC), \underline{\Sym^{k-2}})\otimes \CC_p \right)_{K_p, k}, \]
then we deduce that $\TT'^\wedge_\aux$ is isomorphic to the completed Hecke algebra for (1). By Lemma \ref{lem:CompBaseChange}, $\TT'^\wedge_\aux$ is also the completed Hecke algebra for $\TT'$ acting on 
\[ \left( H^1(Y_{K_p K^p}(\CC), \underline{\Sym^{k-2}})\otimes \QQ_p \right)_{K_p, k}. \]

As in \cite{emerton:local-global-conjecture, emerton:local-global-proof}\footnote{Note that our normalizations for actions are different, so that we obtain a $\Sym^{k-2}$ in the source of the $\Hom$ whereas in \cite{emerton:local-global-conjecture, emerton:local-global-proof} there is a $(\Sym^{k-2})^*$ (cf. \cite[last paragraph of Section~2]{ emerton:local-global-proof}.)},
\[ H^1(Y_{K_p K^p}(\CC), \underline{\Sym^{k-2}})\otimes \QQ_p  \hookrightarrow \Hom_{K_p} (\Sym^{k-2} \QQ_p^2, \wh{H}^1_{K^p}). \]
Thus, if we fix for each $k \geq 2$ a non-zero vector in $(\Sym^{k-2} \QQ_p^2)$, then pairing with these vectors gives Hecke-equivariant injections 
\[ H^1(Y_{K_p K^p}(\CC), \underline{\Sym^{k-2}}) \otimes \QQ_p  \hookrightarrow \widehat{H}^1_{K^p} \]
whose joint image is dense (indeed, it is already so if we fix $k=2$), and thus we deduce from Lemma \ref{lem:CompDense} that $\TT'^\wedge_\aux$ is isomorphic to the completed Hecke algebra for (4), establishing (1)=(4). Then, running the same argument using only weight two modular forms, we find (3)=(4). 

Arguing similarly and using the density of $\GL_2(\ZZ_p)$-algebraic vectors in $\widehat{H}^1_{K^p}$ as established in \cite[Remark 5.4.2]{emerton:local-global-proof} (specifically of the ones which transform locally as $\Sym^{k-2}\QQ_p^2$ for some $k$; we do not need to also allow for arbitrary twists by a determinant), 
we obtain (2)=(4).  
\end{proof}

\begin{remark}\label{remark.final-caraiani-scholze}
From our perspective, instead of $p$-adic modular forms it is perhaps more natural to consider the larger space of $p$-adic automorphic forms given by functions on the Caraiani-Scholze Igusa formal scheme over the ordinary locus, which parameterizes isomorphisms
\[ E[p^\infty] \xrightarrow{\sim} \mu_{p^\infty} \times \bbQ_p/\bbZ_p.
\]
Indeed, an argument nearly identical to the one given in this section for the supersingular Igusa variety but starting with the point 
\[ x=[0:1]\in \bbP^1(\bbQ_p) \]
shows that the completed Hecke algebra of this space of $p$-adic automorphic forms is the same as the one appearing in Theorem \ref{theorem:list-of-isomorphisms}.

Moreover, this space of $p$-adic automorphic forms admits an action of a very large unipotent group at $p$, and using this action one can produce a $\GL_2(\bbA_f^{(p)})$-equivariant projection operator (a type of Kirillov functor) to Katz $p$-adic modular forms which can in turn be used to deduce an isomorphism of completed Hecke algebras -- this will be explained in future work (cf. also \cite{howe:unipotent-circle}). In this way one will be able to obtain a proof of Theorem \ref{mainthm:HeckeIso} that does not pass through singular/\'{e}tale cohomology to show that level and weight families of classical modular forms give rise to the same completed Hecke algebra. 
\end{remark}

\bibliographystyle{plain}
\bibliography{refs}

\end{document}